\documentclass{tac}
\usepackage{xy}
\xyoption{all}

\def\mathcaldef#1{\expandafter\def\csname#1\endcsname{{\cal#1}}}

\def\bb{|\!|}

\def\qq{\quad\quad}
\def\qv{\qq ;\qq}
\def\iso{\,\cong\,}
\def\equ{\,\simeq\,}
\def\inc{\hookrightarrow}
\def\imp{\Rightarrow}
\def\down{\downarrow\!\!}
\def\up{\uparrow\!\!}
\def\la{\langle}
\def\ra{\rangle}
\def\adj{\dashv}
\def\bs{\backslash}
\def\op{^{\rm op}}
\def\st{^*}
\def\ex{\exists}
\def\comp{\pi_0}
\def\p1{$\pi_1$}
\def\ov{\overline}
\def\ul{\underline}
\def\nh{\widehat}
\def\eps{\varepsilon}
\def\tm{\times}
\def\otm{{\,}\otimes}
\def\ootm{{\,}\bar\otimes}
\def\lh{local homeomorphism }
\def\lhs{local homeomorphisms }
\def\pm{perfect map }
\def\pms{perfect maps }

\def\dfs{discrete fibrations }

\def\dofs{discrete opfibrations }
\def\fs{factorization system }
\def\fss{factorization systems }
\def\bfc{balanced factorization category }

\def\ng{neighborhood }
\def\ngs{neighborhoods }
\def\rs{reciprocal stability }
\def\cf{comprehensive factorization }

\def\fcs{finite coverings }

\mathcaldef{C}
\mathcaldef{E}
\mathcaldef{M}
\mathcaldef{S}
\mathcaldef{P}
\mathcaldef{K}
\mathcaldef{B}
\mathcaldef{T}
\mathcaldef{O}
\mathcaldef{D}
\mathcaldef{A}

\mathbfdef{Set}
\mathbfdef{Cat}
\mathbfdef{Gph}
\mathbfdef{Pos}
\mathbfdef{Top}
\mathbfdef{Bim}
\mathbfdef{e}

\mathrmdef{hom}
\mathrmdef{id}
\mathrmdef{ten}
\mathrmdef{Iso}
\mathrmdef{Ar}
\mathrmdef{true}
\mathrmdef{false}
\mathrmdef{lim}
\mathrmdef{colim}
\mathrmdef{exp}

\def\PrX{\Set^{X\op}}
\def\CatX{\Cat\!/\! X}

\def\TopX{\Top\!/\! X}

\def\CX{\C\!/\! X}
\def\CY{\C\!/ Y}
\def\CZ{\C\!/Z}
\def\TX{\T\!/\! X}

\def\BX{\B\!/\! X}
\def\BY{\B\!/Y}
\def\MX{\M\!/\! X}
\def\MY{\M\!/Y}
\def\MP{\M\!/\!P}
\def\MoX{\M'\!\!/\! X}
\def\MoY{\M'\!\!/Y}
\def\MoP{\M'\!/\!P}
\def\Mu{\M/1}
\def\Mou{\M'\!/1}
\def\Eu{\E/1}

\def\EM{{\cal(E,M)}}
\def\EMo{{\cal(E',M')}}
\def\TT{\T}
\def\So{\S_0}
\def\iom{\underline\hom}

\newtheorem{prop}{Proposition}
\newtheorem{corol}{Corollary}
\let\pf\proof
\let\epf\endproof
\def\eq{\begin{equation}}
\def\eeq{\end{equation}}

\author{Claudio Pisani}
\address{via Gioberti 86,\\ 10128 Torino, Italy.}
\title{Balanced category theory II}
\copyrightyear{2009}
\eaddress{pisclau@yahoo.it}

\begin{document}

\maketitle

\begin{abstract}
In the first part, we further advance the study of category theory in a strong \bfc $\C$~\cite{pis},
a finitely complete category endowed with two reciprocally stable \fss such that $\Mu = \Mou$.
In particular some aspects related to ``internal" (co)limits and to Cauchy completeness are considered.

In the second part, we maintain that also some aspects of topology can be effectively synthesized 
in a (weak) \bfc $\TT$, whose objects should be considered as possibly ``infinitesimal" and suitably ``regular"
topological spaces.
While in $\C$ the classes $\M$ and $\M'$ play the role of discrete fibrations and opfibrations, 
in $\TT$ they play the role of \lhs and perfect maps, so that $\Mu$ and $\Mou$ are the 
subcategories of discrete and compact spaces respectively.

One so gets a direct abstract link between the subjects, with mutual benefits.
For example, the slice projection $X/x\to X$ and the coslice projection $x\bs X\to X$, 
obtained as the second factors of $x:1\to X$ according to $\EM$ and $\EMo$ in $\C$, 
correspond in $\TT$ to the ``infinitesimal" neighborhood of $x\in X$ and to the closure of $x$. 
Furthermore, the open-closed complementation (generalized to reciprocal stability) 
becomes the key tool to internally treat, in a coherent way, some categorical concepts (such as (co)limits
of presheaves) which are classically related by duality.

\tableofcontents

\end{abstract}

\section{Introduction}
\label{intro}

In~\cite{pis} we argued that a good deal of basic category theory can be carried out in any 
strong ``balanced factorization category" (bfc). 
Recall that a finitely complete category $\C$ is a bfc if it is endowed with two factorization systems $\EM$ and $\EMo$
which are reciprocally stable:
the pullback of a map in $\E$ (resp. $\E'$) along a map in $\M'$ (resp. $\M$) is itself in $\E$ (resp. $\E'$).
We say that $\C$ is a ``strong" bfc if, furthermore, $\Mu=\Mou$ (the category $\S$ of ``internal sets"). 
We refer to ``weak" bfc's when we wish to emphasize that this condition is not required to hold.
The motivating example of a strong bfc is $\Cat$, with the comprehensive factorization systems: $\M$ and $\M'$ are the
classes of discrete fibrations and opfibrations, while $\E$ and $\E'$ are the classes of final and initial functors,
so that $\Mu=\Mou\equ\Set$ (while $\E/1 = \E'/1$ are the connected categories).

In the first part of the present paper, we review and further develop some aspects of balanced category theory.
In particular, we consider the bifunctors $\ootm_X:\CX\tm\CX\to\S$ and their
restrictions $\otm_X:\MoX\tm\MX\to\S$, where $n\otm_X m := \comp (n\tm_X m)$ is the internal set of components 
of (the total of) the product over $X$ (and reduces in $\Cat$ to the tensor product of the corresponding set-functors).
 
Now, while the bifibrations associated to the factorization systems of the bfc $\C$ 
are summarized, in terms of indexed categories, by the adjunctions
\[
f_!\adj f\st:\CY\to\CX 
\]
\[
\ex_f\adj f\st:\MY\to\MX \qv \ex'_f\adj f\st:\MoY\to\MoX
\]
for any $f:X\to Y$ in $\C$, and the reflections
\[
\down_X\adj i_X:\MX\to\CX \qv \up_X\adj i'_X:\MoX\to\CX
\]
for any $X\in\C$ (and in particular $\comp:=\,\,\down_1 =\,\,\up_1\adj i:\S\to\C$), 
the \rs axiom allows us to obtain also the following ``coadjunction" laws:
\[
f\st n \otm_X m \iso n\otm_Y \ex_f m  \qv n\otm_X f\st m \iso \ex'_f n\otm_Y m    
\]
\[
n\,\ootm_X\, q \iso n\,\otm_X \down_X q \qv  p\,\ootm_X\, m \iso\, \up_X p\otm_X m 
\]
natural in $m\in\MX$ (or $\MY$), $n\in\MoY$ (or $\MoX$) and $p,q\in\CX$.

With this toolkit, we are in a position to straghtforwardly prove familiar properties of
colimits of ``internal-set-valued" maps $m\in\MX$ or $n\in\MoX$, and also that, for any $x:1\to X$ in $\C$, 
there is a bicartesian arrow $\up_X x \to\, \down_X x$ of the bimodule $\ten_X:(\MoX)\op\to\MX$, 
obtained by composing $\otm_X$ with the points functor $\S(1,-):\S\to\Set$. 
Thus the subcategories $\ov X$ of ``slices (projections)" $\down_X x : X/x\to X$ in $\MX$ and $\ov X'$ 
of ``coslices (projections)" $\up_X x : x\bs X\to X$ in $\MoX$ are dual.
Furthermore, under a ``Nullstellensatz" hypothesis, we prove that the bicartesian arrows of $\ten_X$
correspond to the retracts of slices in $\MX$ (or coslices in $\MoX$), so offering an alternative
perspective on Cauchy completion also in the classical case $\C=\Cat$.
It is also shown how these retracts may arise as reflections of figures $P\to X$ whose shape $P$
is an ``atom" (such as the monoid with an idempotent non-identity arrow for $\C=\Cat$). 

In the second part, most of which can be read indipendently from the first one,
we sketch how some relevant aspects of topology can be developed in a bfc too.
While \pms are known to form the second factor of a \fs on the category $\Top$ of topological spaces,
we intend to show that, by replacing $\Top$ with a suitable category $\TT$,
it is reasonable to assume that the same is true for \lhs and that \rs holds therein.

The existence of a reflection $\comp:\T\to\Mu$ in ``sets" suggests that the spaces $X\in\TT$ are ``locally connected",
and in fact the neighborhoods $X/x$ are connected that is, the map $!_{X/x}:X/x\to 1$ is in $\E$.
Some homotopical properties of spaces can be studied through``finite coverings" that is, maps in $\B = \M\cap\M'$;
for instance, a space is ``simply connected" if $!_X^*:\B/1\to\B/X$ is an equivalence.
By the reciprocal stability law, spaces in $\TT$ are also locally simply connected,
so that finite coverings are in fact locally trivial (Corollary~\ref{p.hom}).

Thus we maintain that (weak) bfc's form a common kernel shared by category theory and topology,
and that both the subjects are enlighted by this point of view.
For example, the \rs law allows us, on the topological side, to extend (via exponentiation)
the classical complementarity between open and closed parts to \lhs and perfect maps in $\TT$,
with evident conceptual advantages; 
on the other side, it provides a sort of internal duality for categorical concepts (as sketched above)
which often turns out to be more effective than an ``obvious" duality functor.

\subsection{Outline}

After three preliminary ``technical" sections on bimodules, factorization systems and balanced
factorization categories, and after recalling some concepts of balanced category theory, 
we emphasize in sections~\ref{int}, \ref{ten} and~\ref{sli} 
the central role of the \rs law in treating ``internal aspects" of (balanced) category theory. 
Namely, we study (co)limits of internal presheaves in $\MX$ or $\MoX$, and the role of the retracts
of the representable ones (that is, (co)slice projections).
In the last three sections we sketch the idea of balanced topology; 
in particular, we present some ``evidences" of the fact that the \rs law should hold in an appropriate 
``topological" category $\TT$, in which \lhs and \pms are assumed as the basic concepts.

\section{Bicartesian arrows of bimodules}
\label{bim}

In this section we collect some basic facts about bimodules that will be used in the sequel;
while most of them are well known, others (Proposition~\ref{p.bim}) are new to our knowledge.
We assume that the reader is familiar with the definition of fibration.

Recall that a bimodule $t:X\to Y$ can be seen as a bifunctor $t:X\op\tm Y \to \Set$ or as 
a functor $t:T\to 2$, where $2$ is the arrow category $<\,:0\to 1$.
We pass from one representation to the other, depending on the convenience.

A bimodule is ``representable on the right" if it is a prefibration (or, equivalently, a fibration):
for any $y\in Y$, $t(-,y):X\op\to\Set$ is representable: $t(-,y)\iso X(-,\ov y)$.
Dually, a bimodule is ``representable on the left" if it is an op(pre)fibration.
It is a bifibration iff it is birepresentable, that is corresponds (up to choice) to an adjunction $X\rightharpoonup Y$.

Given a bimodule $t:T\to 2$, if $ga=bf$ is a square in $T$ as below, we write $g (a,b\,) f$.
\[
\xymatrix@R=2pc@C=3pc{
y \ar[d]_g   & x \ar[d]^f \ar[l]_a  \\
y'            & x' \ar[l]_b    \\
1             & 0   \ar[l]              }
\]
If $a$ is opcartesian, then the relation $(\,a,b\,)$ is a function $X(x,x')\to Y(y,y')$
(which in the case of representable bimodules, if $b$ is opcartesian too, becomes the hom-set mapping of a 
corresponding functor $X\to Y$).
In order to graphically emphasize this, when $a$ is opcartesian we write $\la\, a,b\,)$ in place of $(\,a,b\,)$,
and similarly $(\,a,b\,\ra$ if $b$ is cartesian.
If both conditions hold, we have a bijection $\la\, a,b\,\ra$, which in the case of representable bimodules becomes 
the hom-set bijection of a corresponding adjunction $X\rightharpoonup Y$ 
(note that, in that case, the naturality of the bijection is given simply by composition-juxtaposition
of squares). 

On the other hand, if $a=b$ we write $(a)$ in place of $(a,a)$. 
So $g(a)f$ means that $f$ and $g$ are endomorphisms and $ga=af$;
in particular, for identities, $y(a)x$ simply means that $a:x\to y$.
For representable bimodules $y\la a) x$ says that $y$ is the image of $x$ according to a corresponding functor, 
with $a$ as the universal element.

We will be here mainly concerned with the ternary relation $g\la a\ra f$
(or in particular, for identities-objects, $y\la a\ra x$), 
saying that $a$ is bicartesian (or ``biuniversal") and $f$ and $g$ are related (as above) by it.
In that case, we say that say that $f$ and $g$ are ``conjugate" by $a$.
Often we are interested to existentially quantify this relation over some of the three variables;
for example, we write $y\la - \ra x$ if $x$ and $y$ are conjugate by some arrow,
or we say that $x$ is (or has a) conjugate if this holds for some $y\in Y$. 
\begin{remark} [Fixed categories]
\label{r.bim}
It is easily seen by the above remarks that, for any bimodule $t$, there is an equivalence between the full
subcategories $X_t$ and $Y_t$ of conjugate objects in $X$ and $Y$ respectively;
if the bimodule is birepresentable, we get the classical fact that an adjunction
restricts to an equivalence among the full subcategories of objects with isomorphic
units or counits respectively.
(Indeed, the units and counits are conjugate to isomorphisms in $X$ or $Y$, 
and the latter are the bicartesian arrows over $0$ or $1$.)
Similarly, if in the above situation all cartesian arrows are also opcartesian then $Y_t=Y$ 
and so the right adjoint is fully faithful.
\end{remark}
Now we prove that by splitting conjugate idempotents, one gets conjugate objects;
for clarity of notations, we now consider a bimodule $t:X\to X'$ and use primes to denote objects or arrows in $X'$. 
\begin{theorem}   \label{p.bim}
If $e'\la - \ra e$ are conjugate idempotents which split through $y'$ and $y$ respectively,
then $y'\la - \ra y$. More precisely, if $e'\la u \ra e$, $e'=i'r'$ and $e=ir$, then $y'\la r'ui \ra y$.
\end{theorem}
Let us show that $y'\la r'ui\,)y$ that is, that $r'ui$ is opcartesian;
that it is cartesian as well is proved dually.
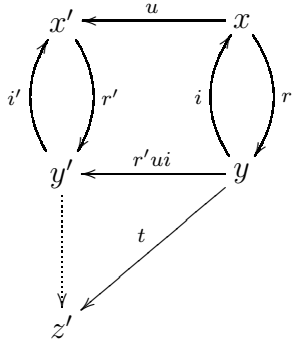
\[
\xymatrix@R=3.5pc@C=4.5pc{
x' \ar@/^1pc/[d]^{r'}  & x \ar@/^1pc/[d]^r \ar[l]_u  \\
y' \ar@{..>}[d]\ar@/^1pc/[u]^{i'}   & y \ar[l]_{r'ui}\ar[ld]_t\ar@/^1pc/[u]^i  \\
z'              &                          }
\]
We need the following
\begin{lemma}   
The retraction $-\circ r',-\circ i':X'(x',z')\to X'(y',z')$ rescrits to a retraction
between the set $F$ of arrows $f':x'\to z'$ such that $f'u = tr$ and the set $G$ of 
arrows $g':x'\to z'$ such that $g'(r'ui) = t$.
\end{lemma}
The theorem follows from the lemma: since $F$ is a terminal set by the hypothesis, the same holds
for its retract $G$, showing that $r'ui$ is cartesian.
\pf
If $f'\in F$ then $f'i'\in G$:

\( f'i'(r'ui) = f'e'ui = f'uei =f'ui = tri = t \)

If $g'\in G$ then $g'r'\in F$:

\( g'r'u = g'r'e'u = g'r'ue = g'(r'ui)r = tr \)
\epf 
\begin{corol} 
Given a bimodule $t:X\to Y$,  
if $X$ and $Y$ are Cauchy complete the same holds for the fixed categories $X_t$ and $Y_t$.
\epf
\end{corol}
In Section~\ref{ten}, we will treat bimodules $X\op\to Y$, that is bifunctors $X\tm Y\to\Set$.
We leave to the reader the simple task of rephrasing the above results to fit this situation.

\section{Factorization systems}
\label{fs}

We assume that the reader is familiar with the basic facts about orthogonality and factorization systems. 
We begin by presenting some properties that will be useful in the sequel and conclude 
by recalling the bifibration associated to a \fs on a finitely complete category.

\begin{prop}  \label{p1.fs}
If $L\adj R:\C\to\C'$ is an adjunction, $Lf\perp g$ iff $f\perp Rg$.
\end{prop}
\pf
By duality, it is sufficient to prove one direction; suppose that $Lf\perp g$
and that the right hand square below commutes.
\[  
\xymatrix@R=4pc@C=4pc{
LA' \ar[r]^{h^*}\ar[d]_{Lf}      & A \ar[d]^g  \\
LB' \ar[r]^{k^*}\ar@{..>}[ur]^u  & B            }
\qq\qq 
\xymatrix@R=4pc@C=4pc{
A' \ar[r]^h\ar[d]_f            & RA \ar[d]^{Rg}  \\
B' \ar[r]^k\ar@{..>}[ur]^{u_*} & RB               }
\]
Then, by the naturality of the transpose bijections, the left square commutes as well, 
giving a unique diagonal $u$; its transpose $u_*$, again by naturality, is easily checked to be 
the desired unique diagonal.
\epf
\begin{prop}   \label{p2.fs}
Let $\EM$ be a \fs on a category $\C$. 
The following are equivalent for a map $e:P\to X$:
\begin{enumerate}

\item
there exists $n:X\to Y$ in $\M$ such that any square $n\circ e =  m\circ l$, with $m\in\M$, has a unique diagonal;
\item
for any triangle $e = m\circ l$, with $m\in\M$, there is a unique section of $m$ extending $l$;
\item
$e\in\E$.
\end{enumerate}
\end{prop}
\pf
The above conditions say, respectively, that the squares below (with $m\in\M$ and the map $n\in\M$ 
in the first one being fixed) have a unique diagonal:
\[  
\xymatrix@R=4pc@C=4pc{
P \ar[r]^l\ar[d]_e         & A \ar[d]^m  \\
X \ar[r]^n\ar@{..>}[ur]^u  & Y            }
\qq\qq
\xymatrix@R=4pc@C=4pc{
P \ar[r]^l\ar[d]_e            & A \ar[d]^m  \\
X \ar@2{-}[r]\ar@{..>}[ur]^u  & X            }
\qq\qq
\xymatrix@R=4pc@C=4pc{
P \ar[r]^l\ar[d]_e         & A \ar[d]^m  \\
X \ar[r]^h\ar@{..>}[ur]^u  & Y            }
\]
(1) $\imp$ (2).
To find the unique section $u$, we compose out with $n$ finding an unique $u$ such that $u\circ e = l$ 
and $n\circ m\circ u = n$: 
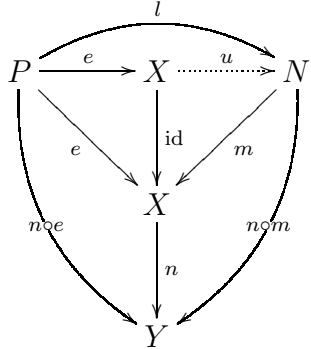
\[
\xymatrix@R=3pc@C=3pc{ 
P \ar[r]^e \ar@/^1.5pc/[rr]^l \ar[dr]_e \ar@/_1.5pc/[ddr]|{n\circ e} &  X \ar[d]^\id \ar@{..>}[r]^u  & N \ar[dl]^m \ar@/^1.5pc/[ddl]|{n\circ m}  \\
                                                                        &  X \ar[d]^n          &                    \\
                                                                        &  Y                   &                     }
\] 
It remains to show that $m\circ u = \id_X$;
this follows from the unicity of the diagonal $m\circ s$ in the square below, 
since $m\circ u\circ e = m\circ l = e$:
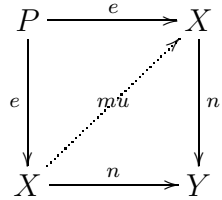
\[  
\xymatrix@R=4pc@C=4pc{
P \ar[r]^e\ar[d]_e           & X \ar[d]^n  \\
X \ar[r]^n\ar@{..>}[ur]|{mu} & Y            } 
\]
(2) $\imp$ (3). 
Consider the factorization $e=me'$ and the uniquely induced diagonal on the left:
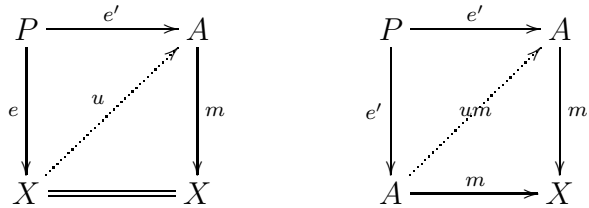
\[
\xymatrix@R=4pc@C=4pc{
P \ar[r]^{e'}\ar[d]_e         & A \ar[d]^m  \\
X \ar@2{-}[r]\ar@{..>}[ur]^u  & X            }
\qq\qq
\xymatrix@R=4pc@C=4pc{
P \ar[r]^{e'}\ar[d]_{e'}     & A \ar[d]^m  \\
A \ar[r]^m\ar@{..>}[ur]|{um} & X            } 
\]
The square on the right shows that also $um=\id$, so that $m$ is an isomorphism.
 
Trivially (3) $\imp$ (1), and the proof is complete. 
\epf

\begin{corol}    \label{p3.fs}
An $\EM$-factorization of a map $f:P\to Y$ in $\C$
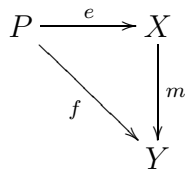
\[
\xymatrix@R=3pc@C=3pc{
P \ar[r]^e\ar[dr]_f  & X \ar[d]^m \\
                     & Y           }
\]
gives both a reflection of $f\in\CY$ in $\MY$ (with $e$ as reflection map) and a coreflection
of $f\in P\bs\C$ in $P\bs\E$ (with $m$ as coreflection map).
Conversely, any such (co)reflection map gives an $\EM$-factorization.
\end{corol}
\pf 
One direction is straightforward. For the converse, note that to say that $e:(P,ne)\to (X,n)$ 
is a reflection of $ne$ in $\MY$ is exactly condition (1) of Proposition~\ref{p2.fs}.
The rest follows by duality.
\epf
\begin{prop}   \label{p4.fs}
Let $\M$ be a pullback-stable class of maps in a finitely complete category $\C$. 
The following are equivalent for a map $e:P\to X$:
\begin{enumerate}
\item
for any triangle $e = m\circ l$, with $m\in\M$, there is a unique section of $m$ extending $l$;
\item
$e\perp m$, for any $m\in\M$;
\item
the pullback functor $e\st:\CX\to \C/\!P$ gives a bijection $\CX(1_X,m)\iso\C/\!P(1_P,e\st m)$, for any $m\in\MX$,
between sections of $m$ and sections of $e\st m$.
\end{enumerate}
\end{prop}
\pf
(1) $\iff$ (2). One implication is trivial. 
For the other one, recall the adjunction $h_!\adj h\st:\CY\to\CX$ and, given $g:B\to Y$, 
denote the corresponding map to the terminal in $\CY$ by $\hat g:(B,g)\to (Y,\id_Y)$. 
Squares in $\C$ with edges $f$, $g$ and $h$ and their diagonals correspond to squares in $\CY$
with edges $h_!\hat f$, $\hat g$ and $\hat h$ and their diagonals:
\[
\xymatrix@R=4pc@C=4pc{
A \ar[r]^l\ar[d]_f         & B \ar[d]^g  \\
X \ar[r]^h\ar@{..>}[ur]^u  & Y            }
\qq\qq
\xymatrix@R=4pc@C=4pc{
(A,hf) \ar[r]^{g_!\hat l}\ar[d]_{h_!\hat f}      & (B,g) \ar[d]^{\hat g}  \\
(X,h)  \ar[r]^{\hat h}\ar@{..>}[ur]^u            & (Y,\id_Y)         }
\]
so that $f \perp g$ in $\C$ iff $h_!\hat f\perp\hat g$ in $\CY$ for any $h:X\to Y$.
By Proposition~\ref{p1.fs}, $h_!\hat f\perp\hat g$ in $\CY$ iff $\hat f\perp h\st\hat g$ $(= \widehat{h\st g})$ 
in $\CX$, that is iff any square as the right hand below has a unique diagonal:
\[
\xymatrix@R=4pc@C=4pc{
(A,f) \ar[r]\ar[d]_{\hat f}            & (h\st B,h\st g) \ar[d]^{h\st\hat g}  \\
(X,\id_X)\ar@2{-}[r]\ar@{..>}[ur]^u         & (X,\id_X)                      }
\qq\qq
\xymatrix@R=4pc@C=4pc{
A \ar[r]\ar[d]_f        & h\st B \ar[d]^{h\st g} \\
X \ar@2{-}[r]\ar@{..>}[ur]^u & X                      }
\]
Since $\M$ is pullback-stable, by the hypothesis the last condition holds for $f=e$ and $g=m$ (for any $m\in\M$) 
so that $e\perp m$, for any $m\in\M$.

(1) $\iff$ (3).
Again by the adjunction $e_!\adj e\st:\CX\to\C/P$, condition (3) says that there is a 
bijection $\CX(1_X,m)\iso\C/P(e, m)$, which is easily seen to correspond to the one of condition (1). 
\epf
\begin{corol}
If $\C$ is finitely complete, Proposition~\ref{p2.fs} holds true for prefactorization systems as well.
\end{corol}
\pf
The implication (2) $\imp$ (3) follows from Proposition~\ref{p4.fs} above.
\epf
\begin{corol}      \label{p6.fs}
Let $\M$ be a pullback-stable class of maps in a finitely complete category $\C$ and $f:X\to Y$ a map in $\C$.
If $f\st:\MY\to\MX$ is fully faithful, then $f\perp m$, for any $m\in\M$.
\epf
\end{corol}
\begin{corol}   \label{p5.fs}
Let $\EM$ be \fs on a finitely complete category.
A map $f:X\to Y$ is in $\E$ iff the functor $f\st:\MY\to \MX$ gives a bijection between the hom-sets 
$\MY(1_Y,m)\iso\MX(1_X,f\st m)$, for any $m\in\MX$.
In particular, if $f\st:\MY\to\MX$ is fully faithful, then $f\in\E$.
\epf
\end{corol}

\subsection{The bifibration associated to $\EM$}

Let $\EM$ a \fs on a finitely complete category.
By restricting the codomain fibration to the arrows in $\M$ we get a subfibration $\M^\to \to\C$
which is a bifibration: the cartesian arrows are the pullback squares again and the opcartesian arrows are
the squares with the top row in $\E$:
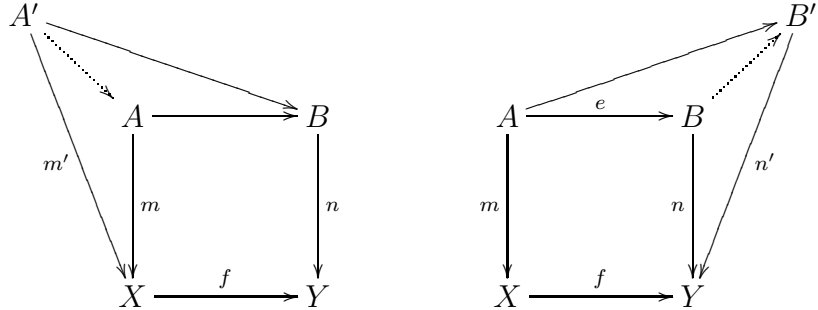
\[
\xymatrix@R=2pc@C=2pc{
A' \ar[dddr]_{m'}\ar@{..>}[rd]\ar[rrrd]         \\
& A \ar[rr]  \ar[dd]^m          &&  B  \ar[dd]^n     \\ \\
& X \ar[rr]^f                   &&  Y                     }
\qq\qq
\xymatrix@R=2pc@C=2pc{
                              &&& B' \ar[dddl]^{n'} \\
A \ar[rr]^e\ar[dd]_m\ar[rrru] &&  B  \ar[dd]_n\ar@{..>}[ru]  \\ \\
X \ar[rr]^f                   &&  Y                     }
\]
(Note that the codomain bifibration $\C^\to \to\C$ itself can be thought of as associated to 
the \fs $(\Iso\C , \Ar\C)$.)

From the indexed point of view, we thus have the family $\MX$, $X\in\C$, of full subcategories $i_X:\MX\inc\CX$,
and adjunctions 
\[ \ex_f\adj f\st:\MY\to\MX \] 
for any arrow $f:X\to Y$ in $\C$.
(No confusion should arise from using the same symbol for both the pullback functor $f\st:\CY\to\CX$
and its ``restriction" $\MY\to\MX$.)
By Remark~\ref{r.bim}, $f\st:\MY\to\MX$ is fully faithful iff any cartesian arrow over $f$ is opcartesian as well 
that is, iff pulling back $f$ along maps in $\M$ one gets maps in $\E$.
Conversely, $\ex_f:\MX\to\MY$ is fully faithful iff squares with the top row in $\E$ are pullbacks.

We also recall that for $p:P\to X$ in $\CX$, $\down_X p := \ex_p 1_P$ is a reflection of $p$ in $\MX$:
\[ \MX(\ex_p 1_P,m) \iso \MP(1_P,p\st m) \iso \C/\!P(1_P,p\st m) \iso \CX(p_!1_P,m) = \CX(p,m) \]
We thus have the adjunction 
\[ \down_X\adj i_X : \MX\to\CX \] 
in which the reflection map (unit) $p\,\to\,\down_X p$ is given by the following opcartesian arrow
with domain $1_P$:
\[
\xymatrix@R=3pc@C=3pc{
P \ar[r]^{e_p}\ar@{=}[d]  &  A  \ar[d]^{\down\,_X p}   \\
P \ar[r]^p                &  X                  }
\]
So, it projects in $\C$ to the first factor $e_p$ of an $\EM$-factorization of $p$, 
while $\down_X p$ is its second factor (see also Corollary~\ref{p3.fs}).
Note also that the bifibration $\MX$ restricts to a ``slices" subopfibration, formed by those objects in $\M^\to$
which admit an opcartesian point. (In $\Cat$, these are the slice projections, so that we obtain
the opfibration corresponding to the ``identity" indexed category; see Section~\ref{cat}.)

\section{Balanced factorization categories}
\label{bfc}

\begin{definition}
A {\bf balanced factorization category} {\rm (bfc)} is a finitely complete category $\C$ with two 
factorization systems $\EM$ and $\EMo$ satisfying the {\bf reciprocal stability law} {\rm (rsl):} 
{\em the pullback of a map in $\E$ (resp. $\E'$) along a map in $\M'$ (resp. $\M$) is itself in $\E$ (resp. $\E'$).}
\end{definition}
(In~\cite{pis}, these were called ``weak" bfc).
If furthermore $\Mu = \Mou$, we say that $\C$ is a {\bf strong} bfc.
\begin{remark}  \label{r.bfc}
Any slice $\CX$ of a bfc is itself a bfc, with the classes $\M_X$, $\M'_X$, $\E_X$ and $\E'_X$
of the maps in $\CX$ mapped to $\M$, $\M'$, $\E$ and $\E'$ by the projection $\CX\to X$;
it is strong iff $X$ is a ``groupoidal object" that is, if $\MX=\MoX$.
\end{remark} 
Typical istances are, for a category $X$, the slice $\CatX$ (see Section~\ref{cat}) and, for a poset $X$, 
the poset $\P X$ of the parts of $X$ with the lower-sets (resp. upper-sets) inclusions as $\M$ (resp. $\M'$).
Both of them are strong if $X$ is a groupoid. 

If $\EM$ is a \fs on a finitely complete $\C$ satisfying the Frobenius law 
that is, maps in $\E$ are pullback-stable along maps in $\M$,
then we obtain a ``symmetrical" bfc by posing $\E' = \E$ and $\M' = \M$; all its objects are
groupoidal and all its slices are symmetrical again. 
An exemple of symmetrical bfc is the category of groupoids, with $\M$ the class of covering maps.
(Other istances of bfc's are presented in~\cite{pis}.) 

We now draw some consequences of the above axioms which will be used in the sequel.
Throughout this section, we assume that $\C$ is a (weak) bfc.
\begin{prop}   \label{p.bfc}
Pulling back an $\EM$-factorization $f = m\circ e$ along a map $n\in\M'$ in $\C$ one gets 
an $\EM$-factorization $n\st f=(m\st n)\st e\circ n\st m$.
\end{prop}
\pf
Consider the pullback squares below. Since $n\st m\in\M$ and $m\st n\in\M'$, the result
follows by applying the rsl to the left one: 
\[
\xymatrix@R=4pc@C=4pc{
A\ar[r]^{e'}\ar[d] & B\ar[d]^{m\st n}\ar[r]^{n\st m}  & C \ar[d]^n  \\
X \ar[r]^e         & Y    \ar[r]^m             & Z          }
\]
\epf
\begin{prop}    \label{p1.bfc}
If $K\in\Mou$ and $e\in\E$ then the map $e\tm K$ is also in $\E$.
\end{prop}
\pf
Considering the pullback squares below, the projection $p$ is in $\M'$ and so by the rsl $e\tm K\in\E$: 
\[
\xymatrix@R=4pc@C=4pc{
X\tm K \ar[r]^{e\tm K}\ar[d] & Y\tm K \ar[d]^p\ar[r]  & K \ar[d]  \\
X \ar[r]^e                   & Y    \ar[r]            & 1          }
\]
\epf
\begin{remark}  
Of course, any property in a bfc (such as the above ones) has a ``dual" property, obtained by exchanging $\M$ with $\M'$
and $\E$ with $\E'$.  
\end{remark}

\begin{prop} [The exponential law]   \label{p2.bfc}
If $m\in\MX$, $n\in\MoX$ and $m^n$ exists in $\CX$, then it is in $\MX$; dually, $n^m\in\MoX$.
\end{prop}
\pf
By Remark~\ref{r.bfc}, we can assume $X=1$: if $S\in\Mu$, $K\in\Mou$ and 
the exponential $S^K$ exists in $\C$, then it is in $\Mu$ that is, $e\perp S^K$ for any $e\in\E$; 
by Proposition~\ref{p1.fs}, this amount to $e\tm K \perp S$ for any $e\in\E$, which follows from Proposition~\ref{p1.bfc}.
\epf
\begin{prop}   \label{p3.bfc}
Suppose that, in the cube below, the bottom, the left and the right faces are pullbacks.
If $e\in\E$ and $n\in\M'$ then $e'\in\E$.
\end{prop}
\[
\xymatrix@R=2pc@C=2pc{
                      &   C \ar[dd]\ar[rr]^{e'}\ar[dl] &               & D \ar[dd]\ar[dl]        \\ 
A \ar[dd]\ar[rr]^>>>>>>e &                           &   B \ar[dd]   &                    \\ 
                      &   Z \ar[rr]\ar[dl]        &               & W \ar[dl]^n        \\ 
X \ar[rr]^f           &                           &   Y           &                      }
\]
\pf
Apply the rsl to the top face, which is a pullback as well.
\epf
\begin{remark}
In~\cite{law70} it was remarked that the Beck and Frobenius conditions do not hold in the eed $\PrX$, $X\in\Cat$ 
(see Section~\ref{cat} below); the above proposition says that the Beck condition does hold 
when restricted to pullback squares, in the base category $\C$, whose right edge is in $\M'$ (and conversely for $\Set^X$).
Thus, we can say that the bifibrations associated to a bfc (see Section~\ref{fs}) satisfy the ``mixed" Beck law.
\end{remark}

\section{$\Cat$ as a strong \bfc}
\label{cat}

Balanced category theory is an abstraction of category theory based on an axiomatization of $\Cat$.
It mainly aims to offer a simple but remarkably powerful conceptual frame in which several 
categorical concepts and properties become quite transparent.
However, it also shows that category theory can be developed, for instance, relatively to a groupoid $X$
that is, in $\CatX$, where the category $\S$ of internal sets (see Section~\ref{int})
is the boolean topos $\PrX \iso \Set^X$ of the coverings of $X$, 
or in the category $\Pos$ of posets (see~\cite{pis}) where $\S = 2$. 

The abstraction is based on the fact that $\Cat$ is a strong bfc with the classes $\M$ of \dfs 
and $\E$ of final functors on one side and the classes $\M'$ of \dofs and $\E'$ of initial functors on the other side.
Recall that $p:P\to X$ is final (resp. initial) iff $\comp (x\bs p) = 1$ (resp. $\comp (p/x) = 1$) for any $x\in X$.
Among final (resp. initial) functors there are the right (left) adjoint ones,
since in this case $x\bs p$ (resp. $p/x$) has an initial (resp. terminal) object.
\begin{remark}  \label{r.cat}
We note that final (resp. initial) functors arise as those which are ``$\A$-asph\'erique" in the sense of~\cite{mal},  
where $\A$ is the ``structures d'asph\'ericit\'e \`a gauche" (resp. ``droite") given by the connected categories;
this fact (which, somewhat surprisingly, is not mentioned there) yields several properties of final and initial functors. 
\end{remark}
The indexed category $\MX \equ \Set^{X\op}$, $X\in\Cat$, was axiomatized (among other things)
in the late sixties by Lawvere as an instance of elementary existential doctrine (eed)
satisfying the ``comprehension scheme". 
So, for example, left Kan extensions appear as existential quantifications left adjoint to substitutions:
$\ex_f\adj f\st: \Set^{Y\op}\to\Set^{X\op}$.

That the bifibration corresponding to this eed is associated to a \fs was shown in~\cite{stw}:

\subsection{The \cf systems}
One easily verifies that $\EM$ and $\EMo$ are the pre\fss generated, respectively, by the codomain  
and the domain functors $t,s:1\to 2$ of the arrow. 
After Section~\ref{fs}, to see that these are in fact \fss it is enough to check 
that $\MX$ is reflective in $\CatX$, which follows by a simple generalization of the Yoneda lemma. 
One also easily checks that the \rs law holds (see~\cite{pis}):
\begin{proposition}
$\Cat$, with the \cf systems, is a strong \bfc in which $\Mu=\Mou\equ\Set$.
\epf
\end{proposition}
\begin{remark}
By Remark~\ref{r.cat}, initial (final) functors are in fact stable with respect to pullbacks along any
(op)fibration (not only the discrete ones); indeed, the latter are smooth functors for any 
asphericity structure \cite{mal}.
Thus, one of the features that distinguishes $\Cat$ among other (strong) bfc's is the fact that
final or initial maps are stable with respect to pullbacks along any projection:
if $e:X\to Y$ is in $\E$ then also $e\tm K:X\tm K\to Y\tm K$ is in $\E$, for any $K\in\Cat$
(and not only for $K\in\Mou$, as in Proposition~\ref{p1.bfc}).
\end{remark}
The following proposition is an example of an effective use of the \rs law; it gives characterizations
of absolutely dense (or ``connected") functors, following~\cite{dense}:
\begin{prop}    
Let $f:X\to Y$ be a functor and let $\,[\alpha]\to Y$ be the interval category of factorizations 
of the arrow $\alpha$ in $Y$ with its projection.
The following are equivalent:
\begin{enumerate}
\item
$f\st[\alpha]$ is connected, for any $\alpha$ in $Y$;
\item
$f$ is locally final: in the pullback square below $e$ is final, for any $y\in Y$;
\eq   \label{d.cat}
\xymatrix@R=3pc@C=3pc{
f/y \ar[r]^e\ar[d] &  Y/y \ar[d]  \\
X \ar[r]^f          &  Y             }
\eeq
\item
$f$ is locally initial: $y\bs f \to y\bs Y$ is in $\E'$ for any $y\in Y$;
\item
$f^*:\MY\to\MX$ is full and faithful;
\item
$f^*:\MoY\to\MoX$ is full and faithful.
\end{enumerate}  
\end{prop}
\pf
First note that $[\alpha]\to Y$ is the composite of a coslice and a slice projection:
\[   
\xymatrix@R=2pc@C=2pc{
[\alpha]\iso \alpha\bs(Y/y) \ar[r] & Y/y \ar[r] & Y }
\]
Thus, in the pullback diagram below $n\in\M'$ and, if $e\in\E$ also $e'\in\E$ by the reciprocal stability law:
\[   
\xymatrix@R=2pc@C=2pc{
f\st[\alpha] \ar[rr]^{e'}\ar[dd]\ar[dr] & & [\alpha] \ar[dd]\ar[dr]^n  \\
                     & f/y \ar[rr]^<<<<<<<<<<e\ar[dl] & &  Y/y \ar[dl]  \\
X \ar[rr]^f                   & & Y                 }
\]
Since $\comp[\alpha] = 1$, also $\comp f\st[\alpha]$ (obtained by factorizing 
$f\st[\alpha]\to 1 = f\st[\alpha]\to[\alpha]\to 1$ according to $\EM$) is $1$.
Conversely, to say that $f\st[\alpha]$ is connected for any $\alpha$ is to say that $e\in\E$, by definition. 
Since condition (1) is self-dual, the equivalence of the first three conditions is proved.

Recalling the adjunction $\ex_f\adj f\st:\MY\to\MX$ of Section~\ref{fs} 
(where $\ex_f m$ is obtained by factorizing $fm$ according to $\EM$, generalizing $\comp:\Cat\to\Set\equ\Mu$), 
Diagram~(\ref{d.cat}) shows that local finality of $f$ is equivalent to the fact that
the counit $\ex_f f\st Y/y\to Y/y$ is an isomorphism, for any $y\in Y$. 
By the properties of the ``Yoneda" inclusion of slices of $Y$ into $\MY$, it is also equivalent to the fact that
the any counit $\ex_f f\st m\to m$ is an isomorphism, that is $f\st:\MY\to\MX$ is full and faithful. 
\epf
\begin{remark} 
After reading Section~\ref{colim} below (see in particular Proposition~\ref{p.colim}), 
it will be clear that Diagram~(\ref{d.cat}) can be interpreted
as exhibiting $y$ as an absolute colimit of $f/y\to Y$, thus explaining the term ``absolute density"; 
see also \cite{pis}, were it is shown that part of the above proposition holds true in any strong bfc.
\end{remark}
\begin{remark} 
If $f:X\to Y$ is also full and faithful, then $f\st:\MY\to\MX$ is an equivalence.
Indeed, in this case the adjoint bimodules corresponding to $f$ are an equivalence in the bicategory $\Bim$
of bimodules, which induces an equivalence between $\Bim(1,X)\equ\Set^X$ and $\Bim(1,Y)\equ\Set^Y$.
Alternatively, recall that the left Kan extension along a fully faithful functor is indeed an extension;
thus a functor $f:X\to Y$ is fully faithful iff the unit $m\to f\st\ex_f m$ is an isomorphism, for any $m\in\MX$,
iff $\ex_f:\MX\to\MY$ is fully faithful.
\end{remark}
\begin{corol} 
Any absolutely dense map is both initial and final.
\end{corol}
\pf
This follows from Corollary~\ref{p5.fs}. Alternatively, note that if $f$ is locally final, 
since $\comp(Y/y) = 1$ also $\comp(f/y) = 1$ (see Diagram~(\ref{d.cat})), that is $f$ is initial.   
\epf
For example, the insertion of a category in its groupoidal reflection is both initial and final.

\section{Slices and colimits in a bfc}
\label{colim}

Throughout this section, we assume that $\C$ is a bfc.
We adopt the general policy of denoting the various concepts in $\C$ as the corresponding ones in $\Cat$.
Thus, for instance, the maps in $\M$ and $\E$ are called discrete fibrations and final maps respectively,
and so on.
 
As in internal category theory, there are two aspects of balanced category theory. 
On the one hand, the objects and arrows of $\C$ are (generalized) categories and functor,
and we can consider concepts such as limits or colimits of maps $f:X\to Y$ and adjunctible maps.
As shown in~\cite{pis} and as we partly recall below, familiar properties (such as the preservation
of limits by adjunctible maps) can be proved therein in a more transparent way.
For these aspects, the \rs law play no real role, so that we could in fact consider this as $\EM$-category theory
(see~\cite{pis2}).

On the other hand (and more interestingly) there are ``internal" aspects, in which objects in $\MX$ or $\MoX$
are considered as (contravariant or covariant) internal-set-valued functors. 
In the next section, we show how the rsl is what makes the internal theory to work.   
(Some ``internal" aspects, however, such as the Yoneda Lemma below, depend only on the \fss axioms.)

\subsection{Slices}

By factorizing an ``object" (point) $x:1\to X$ according to $\EM$ and $\EMo$, we obtain the {\bf slice} and 
the {\bf coslice} projection respectively of $X$ at $x$:
\[
\xymatrix@R=3pc@C=3pc{
1  \ar[r]^{e_x}\ar[dr]_x  & X/x \ar[d]^{\down\,\, x}   \\
                          & X                      }
\qq\qq
\xymatrix@R=3pc@C=3pc{
1  \ar[r]^{e_x'}\ar[dr]_x  & x\bs X \ar[d]^{\up\,\, x}   \\
                           & X                       }
\]
So, as remarked in Section~\ref{fs}, we have $\down_X x = \ex_x 1_1$.
 
One of the consequent universal properties is usually known (in $\Cat$) as the Yoneda Lemma:
\[   
\xymatrix@R=3pc@C=3pc{
1 \ar[r]^{e_x}\ar[dr]_x\ar@/^1.5pc/[rr]^a & X/x \ar[d] \ar@{..>}[r]^u & A \ar[dl]^m \\
                                          & X                         &               }
\]
On the other hand, the slice projection $X/x\to X$ is also the ``biggest" (that is, final) object over $X$ with a 
final point over $x$ (see Corollary~\ref{p3.fs}).

\subsection{Cones and colimits}

Given a map $p:P\to X$ and a point $x$ of $X$, a {\bf cone} $\gamma:p\to x$ (resp. $\gamma:x\to p$) is 
a map in $\CX$ from $p$ to the slice projection $\down_X x$ (resp. coslice projection $\up_X x$):
\[
\xymatrix@R=3pc@C=3pc{
P  \ar[r]^\gamma\ar[dr]_p  & X/x \ar[d]    \\
                           & X              }
\qq\qq                           
\xymatrix@R=3pc@C=3pc{
P  \ar[r]^\gamma\ar[dr]_p  & x\bs X \ar[d]    \\
                           & X              }
\]
A cone $\lambda:p\to x$ (resp. $\lambda:x\to p$) is {\bf colimiting} (resp. {\bf limiting})
if it is universal among cones with domain $p$:
\[   
\xymatrix@R=3pc@C=3pc{
P \ar[r]^\lambda\ar[dr]_p\ar@/^1.5pc/[rr]^\gamma & X/x \ar[d] \ar@{..>}[r]^u & X/y \ar[dl] \\
                                                 & X                         &               }
\]
That is, a colimiting cone gives a reflection of $p\in\CX$ in the full subcategory $\ov X$ 
generated by the slice projections of $X$.
The following property is often taken as a definition of final functors in $\Cat$.
(The converse holds in any $\C$ with ``power objects"; see~\cite{pis2}.)
\begin{prop}  \label{p1.colim}
Precomposing with maps in $\E$ does not affect colimits.
\end{prop}
\pf
If $e:Q\to P$ is in $\E$, then factorizing $p:P\to X$ and $pe:Q\to X$ we get isomorphic factors in $\M$;
thus, $p$ and $pe\in\CX$ have the same reflection in $\MX$ and so also in $\ov X$ (if they exist). 
\epf
The above result can obviously be ``dualized" for limits;
more interestingly, we will show in Section~\ref{int} how the \rs law allows us to {\em internalize} it
(see Proposition~\ref{p.int}).
\begin{remark}
With respect to the classical treatment of (co)limits,
the present approach has several advantages also in the case $\C=\Cat$:
considering the colimit functor on $X$ as a (partial) reflection $\CatX\to X$ makes
the proofs of the following properties quite straightforward (see also~\cite{par}).
\begin{enumerate}
\item
The colimit $x$ of $1_X$, if it exists, is terminal in $X$ (since the reflection $\lambda:1_X\to X/x$ is then an iso);
by Proposition~\ref{p1.colim}, the same is true for any final functor $e:P\to X$. 
\item
The colimit of the empty functor $0\to X$ is an initial object;
if $p:P\to X$ and $q:Q\to X$ have colimits $x_p$ and $x_q$, the colimit of $[p,q]:P+Q\to X$
is $x_p+x_q$.
\item
If $P$ is connected (so that $P\to 1$ is in $\E$) and $p:P\to X$ is constant through $x:1\to X$, then
by Proposition~\ref{p1.colim} $x$ is the colimit of $p$;
similarly, if $p$ is locally constant (that is, factors through $\comp P$) then its colimit is the coproduct 
of the corresponding family.  
\end{enumerate}
\end{remark}
Given a cone $p\to x$ over $X$ and a map $f:X\to Y$, we get a cone $fp\to fx$ by composing
with an opcartesian arrow over $f$ (whose codomain is a slice projection again since it has, 
by composition, an opcartesian point as well):
\[
\xymatrix@R=3pc@C=3pc{
P  \ar[r]^\gamma\ar[dr]_p  & X/x \ar[d] \ar[r]^e & Y/fx \ar[d]  \\
                           & X  \ar[r]^f         & Y            }
\]                          
Thus, we say that $f$ {\bf preserves colimits} if it takes colimit cones $\gamma:p\to x$
to colimiting cones $e\gamma:fp\to fx$.
\begin{proposition} [Absolute colimits; see also~\cite{par}]  \label{p.colim}
If a cone $\gamma:p\to x$ is in $\E$, then it is colimiting and is preserved by any map.
\epf
\end{proposition}

\section{Internal aspects of balanced category theory}
\label{int}

Throughout this section, we assume that $\C$ is a {\em strong} bfc.
Following Section~\ref{fs}, the bfc $\C$ gives rise to two subfibrations of the codomain
bifibration which are themselves bifibrations.
From an indexed (or eed) point of view we thus have adjunctions
\[
f_!\adj f\st:\CY\to\CX 
\]
\[
\ex_f\adj f\st:\MY\to\MX \qv \ex'_f\adj f\st:\MoY\to\MoX
\]
for any $f:X\to Y$ in $\C$.
(No confusion should arise from using the same symbol $f\st$ for three different functors,
since all of them are obtained by pulling back.)

As clearly explained in~\cite{law92}, we thus have varying ``quantities" with both extensive and
intensive aspects. Within the ``gros" categories $\CX$, there are the ``petit" ones of left and right 
``discrete quantities": $i_X:\MX\inc\CX$ and $i'_X:\MoX\inc\CX$. 
(For the ``topological" (weak) bfc $\TT$ of Section~\ref{top}, it would be of course more appropriate
to speak of  ``discrete" and ``compact" quantities or spaces over $X$.)

\subsection{Extensive aspects of discrete quantities}

Since we are now working in a strong bfc, the constant left and right discrete quantities 
coincide: $\S := \Mu = \Mou$; we refer to them as {\bf internal sets}.
Thus we have a {\bf components functor} 
\[ \comp\adj i:\S\to\C \] 
where $\comp X$ is the total (in the sense of $\ex$ or $\ex'$) of the bidiscrete quantity $1_X$,
which can be obtained by factorizing $!_X:X\to 1$ according to $\EM$ or $\EMo$: 
\[ \comp X := \ex_{!_X}1_X = \ex'_{!_X}1_X \]
More generally, we have left adjoints $\comp^X\adj(!_X)\st\circ i$ which can be obtained as
\[ 
\comp^X:=\comp\circ(!_X)_!:\CX\to\S 
\]
that is, if $p:P\to X$, $\comp^X p = \comp P$:
\[
\xymatrix@R=4pc@C=4pc{
P  \ar[r]^e\ar[d]_p  & \comp P \ar[d]  \\
X   \ar[r]^{!_X}     & 1               }
\]
Note that the total $\ex_{!_X}m$ of a left discrete quantity $m:A\to X$ over $X$ 
can be obtained as $\comp A = \comp^X m$ (more precisely, $\comp^X i_X m$, where $i_X:\MX\inc\CX$).
Similarly, $\ex'_{!_X}n = \comp^X n$.

\subsection{Internal-set-valued maps}

Passing now to intensive (that is, contravariant) aspects, for any point $x:1\to X$ we get
the (internal) set $x\st m$ by evaluating a left or right discrete quantity $m$ over $X$.
Furthermore, internal sets are included as constantly varying quantities over $X$ by $!_X\st:\S\to\MX$
(and, of course, evaluating $!_X\st S$ at any $x$ returns $S$ itself). 

Thus a discrete fibration (or opfibration) $m\in\MX$ in $\C$ can be considered as an
``internal-set-valued" map. In this perspective the functor $f\st:\MY\to\MX$
can be seen as precomposition of internal-set-valued maps over $Y$ with $f:X\to Y$,
as is evident from the pullback squares below:
\[
\xymatrix@R=4pc@C=4pc{
x\st f\st m = (fx)\st m \ar[r]\ar[d]  & A \ar[d]_{f\st m}\ar[r]  & B \ar[d]^m  \\
1   \ar[r]^x                          & X \ar[r]^f               & Y           }
\]
In $\Cat$, evaluation of $m$ at $x$ gives of course the value at $x$ of the presheaf associated to $m$.
On the other hand, by the adjunction $\ex_f\adj f\st:\MY\to\MX$, the opcartesian arrow (in $\M^\to\to\C$)
\[
\xymatrix@R=4pc@C=4pc{
A   \ar[r]^e\ar[d]_m  &   B \ar[d]^{\ex_f m}  \\
X   \ar[r]^f          &   Y                   }
\]
corresponds to the left Kan extension along $f$.
In particular, the total of $m:A\to X$ in $\MX$, the internal set $\comp A$, corresponds to the (internal) colimit of $m$.
(Classically, one says that the colimit of a presheaf is given by the components of its category of elements).
Of course, similar considerations hold for discrete opfibrations.

Thus, as noted above, the functor $\comp^X:\CX\to\S$ restricts to give, for discrete fibrations or opfibrations, the
{\bf internal-colimit functors} $\MX\to\S$ and $\MoX\to\S$:
\[
\xymatrix@R=3pc@C=3pc{
\MX \ar[d]_{i_X}\ar[drr]^{\ex_{\,!_X}}    \\
\CX  \ar[rr]^{\comp^X}  &&   \S          \\
\MoX \ar[u]^{i'_X}\ar[urr]_{\ex'_{\,!_X}}   }
\]

\subsection{The role of the \rs law}

Consider the pullback square below with $n\in\M'$,
\[
\xymatrix@R=4pc@C=4pc{
A  \ar[r]^{f'}\ar[d]_{f\st n}  & B \ar[d]^n  \\
X  \ar[r]^f                    & Y           }
\]
Note that $f'$, as a map over $Y$, is the counit $f'=\eps_n:f_!f\st n\to n$ of the adjunction $f_!\adj f\st:\CY\to\CX$.
Thus by applying $\comp^Y:\CY\to\S$ we get a natural transformation $\comp^Y\eps:\comp^Y f_!f\st=\comp^X f\st \to\comp^Y:\MoY\to\S$
(which corresponds ``externally" to the canonical 
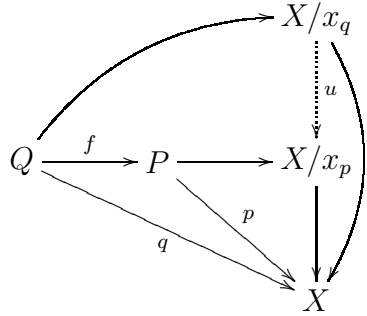
\[   
\xymatrix@R=3pc@C=3pc{
                                      &                   & X/x_q \ar@{..>}[d]^u\ar@/^1.5pc/[dd]     \\
Q \ar[r]^f\ar[drr]_q\ar@/^1.5pc/[urr] & P \ar[dr]^p\ar[r] & X/x_p \ar[d]            \\
                                      &                   & X                      }
\]
where $x_p$ and $x_q$ are colimits of $p$ and $q=pf$ respectively).
Now, if $f\in\E$ then also $f'\in\E$ by the \rs law, and since $\comp:\C\to\S$ takes maps in $\E$ to isomorphisms,
$\comp^Y\eps$ is in fact a natural isomorphism: 
\begin{prop}  \label{p.int}
Precomposition with a final map preserves internal colimits of discrete opfibrations; 
that is, for any $e:X\to Y$ in $\E$, there are isomorphisms 
\[ \comp^X e\st n \iso \comp^Y n \] 
natural in $n\in\MoY$.
\epf
\end{prop}
In particular, the (internal) value of $n\in\MoX$ at a ``final point" $e:1\to X$ in $\E$ 
gives the (internal) colimit of $n$.

Now we apply a similar procedure to obtain other ``coherence" results, supported by the \rs law,
that will be used in the next section.
Considering the diagrams of Section~\ref{bfc}:
\eq   \label{e.int}
\xymatrix@R=4pc@C=4pc{
A\ar[r]^{e'}\ar[d]             & B\ar[d]_{m\st n}\ar[r]^{n\st m}  & C \ar[d]^n  \\
X \ar[r]^e\ar@/_1.5pc/[rr]_p   & Y    \ar[r]^m             & Z          }
\qq
\xymatrix@R=2pc@C=2pc{
                         &   C \ar[dd]\ar[rr]^{e'}\ar[dl] &               & D \ar[dd]\ar[dl]        \\ 
A \ar[dd]^m\ar[rr]^>>>>>>e &                                &   B \ar[dd]_>>>>>{\ex_f m}   &                    \\ 
                      &   Z \ar[rr]^>>>>>>{n\st f}\ar[dl]|<<<<<<{f\st n}   &             & W \ar[dl]^n        \\ 
X \ar[rr]^f           &                           &   Y           &                      }
\eeq
we see that the final maps $e'$ are respectively the components 
\[ \eps_{n,p}:n_!n\st p\to n_!n\st \down p \qv 
\eps_{n,m}:f_!(f\st n)_!(f\st n)\st m\to n_!n\st\ex_f m \] 
of natural transformations between functors
\[ \CZ\tm\M'\!\!/Z\to\CZ \qv
\MX\tm\MoY\to\MY \]
By applying $\comp^Z$ and $\comp^Y$ we get natural transformations 
\[ \comp^Z\eps_{n,p}:\comp^Z(n_!n\st p)=\comp^Z(n\tm_Z p) \to \comp^Z(n_!n\st \down p)=\comp^Z(n\tm_Z\down p) \]
\[ \comp^Y\eps_{n,m}:\comp^Y(f_!(f\st n)_!(f\st n)\st m)=\comp^X(f\st n\tm_X m) \to \comp^Y(n_!n\st\ex_f m)=\comp^Y(n\tm_Y\ex_f m) \]            
where we have used the fact that, if $p:P\to X$ and $f:X\to Y$, $\comp^X p = \comp^Y f_!p$
and, if also $q:Q\to X$, then $\comp^X p_!p\st q = \comp^X q_!q\st p = \comp^X(p\tm_Xq)$:
\[  
\xymatrix@R=4pc@C=4pc{
R  \ar[r]^{q\st p}\ar[d]_{p\st q}\ar[dr]|{p\tm_X q}  & Q \ar[d]^q  \\
P   \ar[r]^p                                         & X               }
\]
Since the $e'$ are in $\E$ and $\comp$ takes final maps to isomorphisms, we get:
\begin{prop}  \label{p1.int}
For any $X\in\C$, there are isomorphisms
\[ \comp^X(n\tm_X p) \to \comp^X(n\tm_X\down p) \]
natural in $n\in\MoX$ and $p\in\CX$. 
For any $f:X\to Y$ there are isomorphisms
\[ \comp^X(f\st n\tm_X m) \to \comp^Y(n\tm_Y\ex_f m) \]  
natural in $m\in\MX$ and $n\in\MoY$.
\epf
\end{prop}

\section{The tensor functor and the internal hom}
\label{ten}

Throughout the section, $\C$ is assumed to be a {\em strong} bfc.
For any $X\in\C$ we define the functor 
\[ \ootm_X:=\comp^X\circ\tm_X:\CX\tm\CX\to\S \]
Thus, if $p:P\to X$ and $q:Q\to X$, we have 
\[ p\ootm_X q = \comp^X(p\tm_X q) = \comp(P\tm_X Q) \] 
By restricting $\ootm_X$ to $\MoX\tm\MX$, we obtain the {\bf tensor functor}
\[ \otm_X:=\ootm_X\circ(i'_X\tm i_X):\MoX\tm\MX\to\S \]
By propositions~\ref{p.int} and~\ref{p1.int}, it immediately follows:
\begin{prop}
For any $X\in\C$, there are isomorphisms
\eq  \label{e1.ten}
n\ootm_X p \iso n\otm_X\down p 
\eeq 
natural in $n\in\MoX$ and $p\in\CX$.
For any $f:X\to Y$ there are isomorphisms
\eq  \label{e2.ten}
f\st n \otm_X m \iso n\otm_Y \ex_f m 
\eeq 
natural in $m\in\MX$ and $n\in\MoY$.
Furthermore, if $e:X\to Y$ is a final map, there are isomorphisms
\eq  \label{e3.ten}
e\st n \otm_X 1_X \iso n\otm_Y 1_Y 
\eeq
natural in $n\in\MoY$. 
\epf
\end{prop}
\begin{remark}
\label{r1.ten}
In fact, both equations~(\ref{e1.ten}) and~(\ref{e3.ten}) follow from Equation~(\ref{e2.ten}).
Indeed, since the Frobenius law $f_!(p\tm_X f\st q) \iso f_!p\tm_Y q$ clearly
holds for the \fs $(\Iso\C,\Ar\C)$ (that is for the codomain fibration), we get isomorphisms 
\[ f_!p\ootm_Y q \iso p\ootm_X f\st q \]
natural in $p\in\CX$, $q\in\CY$, for any $f:X\to Y$. So, if $p:P\to X$ and $n\in\MoX$,
\[ n\ootm_X p \iso n\ootm_X p_!1_P \iso p\st n\ootm_P 1_P \iso p\st n\otm_P 1_P \iso n\otm_X \ex_p 1_P \iso n\otm_X\down p \]
If $e:X\to Y$ is in $\E$, then $\down_Y e = 1_Y$ so that 
\[ e\st n \otm_X 1_X \iso n \otm_Y \ex_e 1_X \iso n \otm_Y \down_Y e \iso n\otm_Y 1_Y \]
\end{remark}

\subsection{The coadjunction laws and the tensor-hom duality}

We have so obtained some ``coadjunction" laws which, remarkably, are the exact counterparts of 
the genuine adjunction laws constituting the logic of the bifibrations originated by $\EM$ and $\EMo$.
Let us emphasize this sort of duality:

\eq     \label{e4.ten}
\begin{array}{c}
\MX(m,f\st  m') \\ \hline
\MY(\ex_f m,m') 
\end{array}
\qq
\begin{array}{c}
\MoX(n,f\st n') \\ \hline
\MoY(\ex_f n,n') 
\end{array}
\qv
\begin{array}{c}
f\st n \otm_X m \\ \hline
n\otm_Y \ex_f m  
\end{array}
\qq
\begin{array}{c}
n\otm_X f\st m \\ \hline
\ex'_f n\otm_Y m  
\end{array}
\eeq

\eq     \label{e5.ten}
\begin{array}{c}
\CX(p,m) \\ \hline
\MX(\down_X p,m) 
\end{array}
\qq
\begin{array}{c}
\CX(q,n) \\ \hline
\MoX(\up_X q,n) 
\end{array}
\qv
\begin{array}{c}
p\,\ootm_X\, m \\ \hline
\up_X p\otm_X m  
\end{array}
\qq
\begin{array}{c}
n\,\ootm_X\, q \\ \hline
n\,\otm_X \down_X q 
\end{array}
\eeq
Furthermore, from the~(\ref{e4.ten}) above and ``surjectivity" of final maps ($\ex_e1_X = \down e = 1_Y$), 
it follows that, if $e\in\E$ and $i\in\E'$, 
\eq    \label{e6.ten}
\begin{array}{c}
\MX(1_X,e\st m) \\ \hline
\MY(1_Y,m) 
\end{array}
\qq
\begin{array}{c}
\MoX(1_X,i\st n) \\ \hline
\MoY(1_Y,n) 
\end{array}
\qv
\begin{array}{c}
1_X\otm_X i\st m \\ \hline
1_Y \otm_Y m  
\end{array}
\qq
\begin{array}{c}
e\st n\otm_X 1_X \\ \hline
n \otm_Y 1_Y 
\end{array}
\eeq
(Note that, by Corollary~\ref{p5.fs}, the converse holds for the left ones.) 
Now, as explained before, the internal set $1_X\otm_X m = \comp^X m$ can be seen as the colimit 
of the internal-set-valued map $m\in\MX$; 
``dually", the (external) set of sections $\MX(1_X,m)$ can be seen as the limit of $m$. 
We will discuss in Section~\ref{lim} below under which hypothesis this limit can be internalized too.
Of course, in $\Cat$ both give the usual limit and colimit of (the presheaf corresponding to) $m$.

More generally, $p\ootm_X m$ can be seen as an internal-set-valued way to ``test" the quantity $m$,
``dual" to the standard set-valued testing by figures of shape $p\in\CX$.
As discussed at lenght in~\cite{pis1} and~\cite{pis4}, the tensor functor $\ootm$ can be seen as a sort
of ``meets" predicate, so that $p\ootm_X m$ gives (the internal set of)
the ways in which $p$ meets $m$ (so as $\CX(p,m)$ gives the ways in which $p$ is contained in, 
or belongs to, $m$). 
In the strong bfc $\Pos$ mentioned at the beginning of Section~\ref{cat} (that is the example at the end of~\cite{pis}),
one has $\S=2=\{\true\vdash\false\}$ and $\otm_X:\MoX\tm\MX\to 2$ is indeed the two-valued meets predicate 
for upper and lower subsets of the poset $X$: $n\otm_X m = \true$ iff $n\cap m$ is non-void. 
\begin{remarks}    \label{r2.ten}
\begin{enumerate}
\item
In this perspective, the internal colimit functor 
\[  1_X\otm_X - \,\,=\,\, \comp^X : \MX \to \S  \]
becomes a ``non-void" predicate (and similarly $\MX(1_X,-):\MX\to\Set$ is a ``whole" predicate).
Preservation of internal colimits (of discrete opfibrations) by precomposition with final maps 
(equations~(\ref{e6.ten}) above) then becomes the fact that the ``surjectivity" of $e:X\to Y$ ($\ex_e 1_X = 1_Y$) 
imply that taking inverse images preserves (and reflects) non-voidness.
Similarly, preservation of limits (of discrete fibrations) says that taking inverse images 
(preserves and) reflects wholeness.
\item
Conversely, the ``meets" and the ``belongs to" predicates can be reduced to the ``non-void"
(colimit) and the ``whole" (limit) predicates by the (co)adjunction laws:
\[ p\ootm_X m \iso 1_P\otm_P p\st m  \qv \CX(p,m) \iso \MP(1_P,p\st m) \]
\item
If $x:1\to X$ is a point, $x\ootm_X m \iso x\st m$ is the internal value of $m$ at $x$ as discussed above, 
while $\CX(x,m)$ is the set of points of the total of $m$ which are (in the fiber) over $x$.
\item
From equations~(\ref{e4.ten}) and~(\ref{e5.ten}) we get the classical formulas
for the left Kan extension $\ex_f m$ of $m\in\MX$ along $f:X\to Y$:
\[ 
y\st \ex_f m \iso y\ootm_Y\ex_f m \iso\,\,\up_Y y\otm_Y\ex_f m \iso f\st\up_Y y\otm_X m 
\]
The last term being $\comp^Y(f\st\up_Y y\tm_X m) = \comp^X(m\st(f\st\up_Y y))$, in $\Cat$ we get
the coend or the colimit formula respectively (see the second of diagrams~(\ref{e.int})).
\end{enumerate}
\end{remarks}

\subsection{Internal limits and internal hom}
\label{lim}

We have argued above that there is a sort of duality between the hom and the tensor functors;
on the other hand there is a great difference: while the latter is valued in $\S$, the former
is valued in $\Set$; to better compare them we need either to unenrich $\otm$ to $\Set$ or to enrich
$\hom$ to $\S$. The first option will be followed in Section~\ref{spl}, where we consider the ``$\ten$" 
bimodule, obtained by composing $\otm$ with the points functor;
we now briefly consider the other one.

So as we obtained internal colimits of discrete (op)fibrations by restricting the left adjoints
$\comp^X\adj(!_X)\st\circ i:\S\to\CX$, and then used these to define the (more general) internal
tensor functor, we now need to {\em assume} the right adjoints 
\[ |-|_X:\CX\to\S \]
to the ``constant" inclusions, which restricted to (op)fibrations give the {\bf internal limit} 
(or ``internal sections") functors.
(These functors in fact exist in $\Cat$ where, since $\S=\Set$ and $!_X\st S = S\cdot 1_X$,
they are the sections functors $|\,p\,|_X = \CX(1_X,p)$.)
Assuming furthermore that, as in $\Cat$, discrete fibrations $m\in\MX$ and opfibrations $n\in\MoX$
are exponentiable in $\CX$ (so that also $n\tm_X m$, $m^p$, etc. are exponentiable),
it is natural to define the ``internal hom" (partial) functors (over $X$) 
as the exponential followed by internal sections:
\[ \iom_X(p,q) = |\,q^p\,|_X \in\S \]     
While $\iom_X$ may be not defined on the whole $\CX\tm\CX$, it is of course defined when the second 
component is exponentiable. 
We denote by $\iom_{\MX}:\MX\tm\MX\to\S$ the restriction of $\iom_X$.
(Note that $m^{m'}$, with $m,m'\in\MX$, may be not in $\MX$.)
\begin{prop}
If the strong bfc $\C$ admits internal limits $|-|_X:\CX\to\S$, we have the ``internal adjunctions"
\[ \iom_X(p,m) \iso \iom_{\MX}(\down_X p,m) \qv \iom_{\MX}(m,f\st l) \iso \iom_{\MY}(\ex_f m,l) \]
natural in $p\in\CX$, $m\in\MX$ and $l\in\MY$; the same holds of course for discrete opfibrations.
\end{prop}
\pf
\[ \S(S,|m^p|_X) \iso \CX(!_X\st S, m^p) \iso \CX(!_X\st S\tm p, m) \iso \MX(\down(!_X\st S\tm_X p), m) \iso \]
\[ \iso \MX(!_X\st S\tm_X \down\,p, m) \iso \CX(!_X\st S, m^{\down\,\,p}) \iso \S(S,|m^{\down\,\,p}|_X)       \]
where, in the passage from the first to the second row, we have applied Proposition~\ref{p.bfc} to 
the $\EM$-factorization of $p$, since the constant $!_X\st S$ is a discrete bifibration.
The second deduction is similar:
\[ \S(S,|(f\st l)^m|_X) \iso \CX(!_X\st S,(f\st l)^m) \iso \MX(!_X\st S\tm m,f\st l) \iso  \]
\[ \iso \MY(\ex_f(!_X\st S\tm_X m), l) \iso \MY(\ex_f(f\st!_Y\st S\tm_X m), l) \iso \] 
\[ \iso \MY((!_Y\st S\tm_Y \ex_f m), l) \iso \CY(!_Y\st S,l^{\ex_f m}) \iso \CY(S,|l^{\ex_f m}|_Y) \]
here, in the passage from the second to the third row, we have applied the mixed Beck law 
(see Proposition~\ref{p3.bfc} and the second of diagrams~(\ref{e.int})).
\epf
\begin{corollary}
Final maps preserve internal limits of discrete fibrations.
\end{corollary}
\pf
\( \iom_{\MX}(1_X,e\st m) \iso \iom_{\MY}(\ex_e 1_X,m) \iso \iom_{\MY}(1_Y,m) \)
\epf
Thus, as for internal colimits, in order to coherently internalize limits and hom (``natural transformations")
of internal-set-valued maps we need the rsl in an essential way.
We conclude this section by comparing the $\iom_{\MX}$ and $\iom_{\MoX}$ (and $\iom_X$) with $\otm_X$ (and~$\ootm_X$),
obtained as the horizontal compositions in the diagrams below:
\[
\xymatrix@R=3pc@C=3pc{
&&           \MX \ar[d]_{i_X}\ar[drr]^{\ul{\colim}_X}\ar@/_1.5pc/[dll]_{\la\,1_X,\,\MX\ra}    \\
\MoX\tm\MX\ar[r]^{i'_X\tm\, i_X} & \CX\tm\CX\ar[r]^{\tm_X} & \CX  \ar[rr]^{\comp^X}  &&   \S          \\
&&           \MoX \ar[u]^{i'_X}\ar[urr]_{\ul{\colim}'_X}\ar@/^1.5pc/[ull]^{\la\,\MoX, 1_X\ra}   }
\] 

\[
\xymatrix@R=1pc@C=3pc{
&&           \MX \ar[dd]_{i_X}\ar[ddrr]^{\ul{\lim}_X}\ar@/_1pc/[dll]_{\la\,1_X,\,\MX\ra}    \\
\MX\tm\MX\ar@/_.5pc/[dr]^{i_X\tm\, i_X}                                                      \\
                                  & \CX\tm\CX\ar[r]^{\exp_X} & \CX  \ar[rr]^{|-|_X}  &&   \S          \\
\MoX\tm\MoX\ar@/^.5pc/[ur]_{i'_X\tm\, i'_X}                                                            \\
&&           \MoX \ar[uu]^{i'_X}\ar[uurr]_{\ul{\lim}'_X}\ar@/^1pc/[ull]^{\la\,1_X,\,\MoX\ra}   }
\]

\section{Retracts of slices}
\label{sli}

\subsection{Components and the Nullstellensatz hypothesis}
\label{null}

As pointed out by Lawvere in several papers,
for categories of cohesion $\C$, whose objects are to be thought of as spaces of some kind,
there is a basic chain of adjoints $p_!\adj p\st\adj p_*:\C\to\S$ (with suitable properties), 
contrasting it with a category $\S$ of (relatively) discrete spaces.
In that situation, he refers to the Nullstellensatz condition as the requirement that (assuming $p\st$ fully faithful)
the natural map $p_*X\to p_!X$, from the points functor to the components (or ``pieces") functor,
is an epimorphism. 
In our setting, we have $\comp\adj i:\S\to\C$, but we do not assume in general a further right adjoint.
Notwithstanding, we will use a weak form of the Nullstellensatz:
if we denote by 
\[ \bb -\bb :=\C(1,-):\C\to\Set \] 
the (external) points functor, and by $[-]_X:X\to\comp X$ the unit of the (internal) components reflection,
we require that the mapping 
\[ \bb [-]_X\bb :\bb X\bb \to\bb \comp X\bb \] 
is surjective, for any $X\in\C$.

Note that for any element $s\in \bb \comp X\bb $ (that is, $s:1\to\comp X$) of the set of 
components of $X$, we have a ``component" $\ul{[s]}\inc X$, that is the subobject given by the 
following pullback:
\eq    \label{e7.ten}
\xymatrix@R=3pc@C=3pc{
\ul{[s]} \ar[d]\ar[r]   & X \ar[d]    \\
1 \ar[r]^s              & \comp X      }
\eeq 
Note also that a figure $p:P\to X$ belongs to (that is, factors through) a component,
iff the composite $[p]:P\to\comp X$ is constant.
In particular, any figure with a connected shape that is, with $!_P\in\E$ (for instance a point), 
belongs to a component. Thus, the Nullstellensatz condition 
\[  \bb \comp X\bb  = \{[x]\,|\,x:1\to X\} \]
may be rephrased by saying that each component has a point (which belongs to it). 
Furthermore, for a map $f:X\to Y$, the corresponding mapping $\bb \comp X\bb \to\bb \comp Y\bb $
acts as $[x]\mapsto [fx]$.

\subsection{The bimodule ten}
\label{spl}

It is well known that the rectracts of slices (representable presheaves) in $\MX$ have
an important role in $\Cat$; for instance, they generate the Cauchy completion of $X$
and can be characterized in several ways.
In order to develop a similar analysis in $\C$, we need to consider the ``unenrichment" 
mentioned in Section~\ref{lim}, by taking the points of $\otm$; 
namely we define the bimodules $\ten_X : (\MoX)\op \to \MX$ by 
\[ \ten_X(n,m):= \bb \, n\otm_X m\,\bb  \]
 
\begin{prop}  \label{p1.ten}
For any $x:1\to X$ there is a bicartesian arrow $\up_X x \to \,\down_X x$ for $\ten_X$.
\end{prop}
\pf
Recalling the notations of Section~\ref{bim}, we begin by showing that $\up_X x ( - \ra \down_X x$, 
that is that $\ten_X(\up_X x,-):\MX\to\Set$ is represented by $\down_X x$:
\[ 
\bb  \up_X x \otm m\, \bb  \iso \bb \,\ex'_x1_1 \otm m\, \bb \iso \bb 1_1 \otm x\st m\, \bb  \iso \bb \, x\st m\, \bb  \iso \]
\[ \iso \S(1,x\st m) \iso \MX(\ex_x1_1,m) \iso \MX(\down_X x,m)
\]
Now, it is easy to see that the cartesian arrow (universal element) is given by the component
$[\la e'_x, e_x \ra]$ of the element $\la e'_x, e_x \ra:1\to x\bs X\tm_X X/x$
(with $e_x\in\E$ and $e'_x\in\E'$):
\[
\xymatrix@R=3pc@C=3pc{
1 \ar[dr]\ar@/^1.5pc/[drr]^{e_x}\ar@/_1.5pc/[ddr]_{e'_x}    \\
& x\bs X\tm_X X/x \ar[d]\ar[r]   & X/x \ar[d]^{\down\,\, x}   \\
& x\bs X  \ar[r]^{\up\,\, x}     & X             }
\] 
Thus, by symmetry, the cartesian arrow is also opcartesian, and the proof is complete.
\epf
\begin{corollary}
The full subcategories $\ov X\inc\MX$ and $\ov X'\inc\MoX$, generated by the slices and 
the coslices projections respectively, are dual.
\epf
\end{corollary}
We are now in a position to prove:
\begin{prop}   \label{p2.ten}
Under the Nullstellensatz hypothesis, all the (op)cartesian arrows of $\ten_X$ are in fact bicartesian,
and the conjugate objects in $\MX$ ($\MoX$) are the retracts of (co)slices projections.
So the fixed categories for $\ten_X$ are the Cauchy completions of $\ov X\inc\MX$ and $\ov X'\inc\MoX$.
\end{prop}
\pf
Observe that, by the hypothesis, for $n:D\to X \in\MoX$ and $m:A\to X \in\MX$, 
\[ \ten(n,m) = \{d\otm a \,|\, d:1\to D, a:1\to A, nd = ma \} \]
where we pose, as usual for $\Cat$, $d\otm a:= [\la d,a \ra]$.
Thus, if $n\la s )m$ is cartesian, we have that $s = v\otm w\in\ten(n,m)$ is over an $x=nv=mw\in X$. 
Let $i:m\to\,\down x$ be the unique map in $\MX$ such that $\ten_X(n, i):v\otm w\mapsto v\otm e_x$,
and $r:\,\down x \to m$ be the unique map in $\MX$ such that $r:e_x\mapsto w$ (where $\down x\circ e_x:1\to X$
is an $\EM$-factorization of $x$).
Then, $\ten_X(n, ri):v\otm w\mapsto v\otm w$ and the cartesianess of $v\otm w$ implies $ri=\id_m$
that is, $m$ is a retract of $\down x$.

By Proposition~\ref{p1.ten} above, $\up x\,\,\la - \ra\down x$ and since $\MoX$ is finitely complete,
the idempotent $e'$, conjugate to $e=ir$, splits as $e'=i'r'$.
Thus, by Proposition~\ref{p.bim}, $m$ has a conjugate in $\MoX$ and the result follows.
\epf

\subsection{Atoms}
\label{at}

Intuitively, an object $P\in\C$ is an atom if it is so small that any non-void open or closed part over it
is the whole $P$, and yet so big that the whole $P$ is itself non-void (see also~\cite{pis1} and~\cite{pis4}).

Now (see Remark~\ref{r2.ten}) $\bb \comp^P m \bb = \bb 1_P\otm_P m \bb = \ten_P(1_P,m)$ 
can be seen as the (external) truth value of the 
``non-voidness" of $m\in\MP$, while  $\hom_P(1_P,m)$ is the (external) truth value of its ``wholeness" 
(where for simplicity we denote by $\hom_P$ the hom-functor on $\C/\!P$ or also its restriction to $\MP$ or $\MoP$).
Thus we formalize the above idea by the conditions
\[ \ten_P(1_P,m)\iso \hom_P(1_P,m)  \qv  \ten_P(n,1_P) \iso \hom_P(1_P,n) \]
(for $m\in\MP$ and $n\in\MoP$) which express the fact that the (external) limit and colimit functors, 
for discrete fibrations and opfibrations, are isomorphic.
In fact, the two conditions are equivalent because, for this particular case, the results of 
Proposition~\ref{p2.ten} can be summarized in the following corollary-definition:
\begin{proposition}
Under the Nullstellensatz hypothesis, the following are equivelent for an object $P$ of $\C$:
\begin{enumerate}
\item
$P$ is an {\bf atom};
\item
$ 1_P\la - \ra 1_P $, for $\ten_P$;
\item
$\ten_P(1_P,m)\iso \hom_P(1_P,m)$, naturally in $m\in\MP$;
\item
$\ten_P(n,1_P) \iso \hom_P(1_P,n)$, naturally in $n\in\MoP$.
\epf
\end{enumerate}
\end{proposition}
A typical case of an atom is an object $X\in\C$ with a ``zero" point $x:1\to X$ in $\E\cap\E'$,
since in that case $x\bs X \iso X/x \iso 1_X$.
In particular, the terminal object itself $1\in\C$ is an atom.

For $\C = \Cat$, the above conditions are related again to the Nullstellensatz condition, 
now referred to the (colimit and limit) adjunctions 
\[ p_!\adj p\st\adj p_*:\PrX\equ\MX\to\Set \]
(see~\cite{rey04}).
Indeed, $X$ is an atom iff $p_!\iso p_*$ or, equivalently, if the same holds for $p_!,p_*:\Set^X\to\Set$.
(Note that atoms in $\Cat$ are connected since 
\[ \comp X = \bb \comp X \bb = \ten_X(1_X,1_X) = \hom_X(1_X,1_X) = 1 \]
so that the corresponding $p\st$ is fully faithful.)
The most relevant instance of a (non point-like) atom in $\Cat$ is the monoid $\e$ with an unique 
idempotent non-identity arrow $e$. 
Indeed, the points-sections-limit functor and the components-colimit functor $\Set^\e\to\Set$ are
isomorphic (to the fixed points of the endomapping associated to $e$).
\begin{remark}
Of course, in presence of internal limits $|-|_X:\CX\to\S$, as discussed in Section~\ref{lim}, 
one may define ``internal atoms" that is, objects $P\in\C$ such that $\comp^P m\iso |m|_P$,
naturally in $m\in\MP$ or $m\in\MoP$.
\end{remark}
\begin{proposition}
The reflection $\down_X x$ (resp. $\up_X x$) of figures $x:P\to X$ with atomic shapes
are in the fixed category of the conjugate objects in $\MX$ (resp. $\MoX$). 
\end{proposition}
\pf
If $P\in\C$ is an atom and $x:P\to X$, the adjunction and (unenriched) coadjunction laws give rise to isomorphisms 
\[ \ten_X(\up_X x,m) \iso \ten_X(\ex'_x 1_P,m) \iso \ten_P(1_P,x\st m) \iso \] 
\[ \iso \hom_P(1_P,x\st m) \iso \hom_X(\ex_x 1_P,m) \iso \hom_X(\down_X x,m) \]
natural in $m\in\MX$ (and similarly for $n\in\MoX$). 
Thus the result follows from Proposition~\ref{p2.ten}.
\epf
In $\Cat$ one so gets in fact {\em all} conjugate presheaves: indeed, a retract $m$ of the representable $\down x$ 
can be obtained as the reflection $\down e$ of the atomic figure $e:\e\to X$ which represents the 
corresponding idempotent ($m\iso\,\down e$ since both of them split the same idempotent in $\MX$; see also~\cite{pis0}).
Another property which characterizes conjugate presheaves in $\Cat$ is the cocontinuity of
the functor $\PrX\to\Set$ represented by it:
\[ \PrX(\down e,-) \iso\hom_X(\down e,-) \iso \hom_X(e,-)\iso \ten_X(e,-)\iso \ten_X(\up e,-) \] 
which gives, for $m\in\PrX$, the elements of $mx$ fixed by $me:mx\to mx$.
For a general bfc $\C$, we have a partial internal version of that fact, as a consequence of the following:
\begin{proposition}  [Complements]
If the constant bifibrations $!_X\st S$, with $S\in\S$, are exponentiable in $\CX$, then
the functor $n\otm_X - :\MX\to\S$ has a right adjoint $(!_X\st -)^n$. 
\end{proposition}
\pf
First note that, by the exponential law (Proposition~\ref{p2.bfc}), $(!_X\st -)^n$ is indeed
valued in $\MX$. Then we have:

\( \S(n\otm_X m,S) \iso \S(\comp^X(n\tm_X m),S) \iso \CX(n\tm_X m,!_X\st S) \iso \MX(m,(!_X\st S)^n) \)
\epf
As discussed at lenght in~\cite{pis1},~\cite{pis0} and~\cite{pis4} (see also Section~\ref{bt} below),
this right adjoint well deserves to be called the ``complement" of $n$.
Now, if $\down e$ is a conjugate object and the points functor $\bb - \bb:\S\to\Set$ preserves itself 
colimits, then the same holds for 
\[ \hom_X(\down e,-) \iso \ten_X(\up e,-) \iso \bb \up e \otm_X - \bb \]

\section{Conclusion of the first part}
\label{con1}

We hope to have shown that the adjunction and coadjunction laws associated to a strong \bfc are a
powerful tool to synthetically treat some basic aspects of category theory.

In this ``categorical logic of categories",
a straightforward common generalization of the ``meets" predicate and of the (internal) colimit functor,
namely the tensor functor, is related to the hom functor (generalizing the ``belongs to" predicate 
and the sections or limit functor) by a useful sort of duality, which is disciplined by the \rs law. 
Further suggestions can be drawn by comparing it with the weaker logic associated to ``topological" 
weak bfc's, which is briefly illustrated in the sequel.

\section{Universal properties in topology}
\label{bt}

It is commonly acknowledged that the main reason of the effectiveness of category theory is its role 
as a language apt to define and elaborate the universal properties which pervade mathematics.
For instance, the universal definition of product gives, in $\Set$ and $\Set\op$, the objectified
version of product and sum of natural numbers, and the right adjoint to $X\times-$ gives exponentials
(and in general implies the distributive law).
Shifting from $\Set$ to $\P X$ (the slice $\Set/X$ restricted to monomorphisms), one similarly gets the boolean algebra
of the parts of $X\in\Set$ (with implication as exponential).

Our present aim is to sketch how some of the universal properties that pervade topology can be used to organize 
and guide our topological thinking.

\subsection{Orthogonality in topology}

Let us consider the concepts of connectedness and density.
In the category $\Top$ of topological spaces, a space $X$ is connected iff any map to a discrete space is constant,
that is if the map $X\to 1$ is orthogonal to $S\to 1$, for any discrete $S$:
\[  
\xymatrix@R=4pc@C=4pc{
X  \ar[r]\ar[d]          &  S \ar[d] \\
1  \ar[r]\ar@{..>}[ur]   &  1         }
\]
Shifting from $\Top$ to $\P X$ (the slice $\Top/X$ restricted to monomorphisms), 
and replacing discrete spaces with closed parts, we get density: 
a part $P$ of $X$ is dense iff any map (that is, inclusion in) to a closed part $D$ is constant
(that is, it factors through the terminal part $X\in\P X$):
\[  
\xymatrix@R=4pc@C=4pc{
P  \ar[r]\ar[d]          &  D \ar[d] \\
X  \ar[r]\ar@{..>}[ur]   &  X         }
\]
In $\Top$, \lhs (resp. perfect maps) to a space $X$ can be seen both as variable discrete 
(resp. compact) spaces over $X$ and as generalized (non monomorphic) open (resp. closed) parts of $X$.
Thus, one is led to consider \lhs and \pms as the basic concepts, and to investigate which are the general 
counterparts of the above orthogonality conditions.

\subsection{Factorization systems}

It is known that perfect (that is, proper and separated) maps are the second factor of a 
factorization system $\EMo$ on $\Top$ (which generalizes the Stone-Cech compactification; see e.g.~\cite{clem96}).
Since $\K = \Mou$ is the subcategory of compact (separated) spaces, $\E'/1$ includes the codiscrete spaces. 

On the other hand, local homeomorphisms are not the second factor of a factorization system $\EM$ on $\Top$: 
assuming that they are so corresponds intuitively to assume both some local connectedness
property on spaces, and the existence infinitesimal neighbouring spaces.
Indeed, $\EM$-factorization gives reflections $\CX\to\MX$ and in particular $\comp:\C\to\S$,
where $\S = \Mu$ is the category of ``internal sets" or ``discrete spaces".
In $\Top$, only the (weakly) locally connected spaces (which are the sum of their components)
have such a reflection.
By considering instead the ``opposite" case of monomorphic figures $P\inc X$, one should obtain the smallest open part
including the figure, that is the (``infinitesimal") neighborhood of $P\inc\nh P\inc X$ of $P$ in $X$:
\[   
\xymatrix@R=4pc@C=4pc{
P  \ar[r]\ar[d]             &  O \ar[d] \\
\nh P  \ar[r]\ar@{..>}[ur]  &  X         }
\]
The same diagram shows that, if $O\to X$ is any map in $\M$ and the figure $p:P\to X$ lifts to $q:P\to O$, 
then its ``neighborhood" $\nh p:\nh P\to X$ lifts uniquely to give a neighborhood $\nh q:\nh P\to O$ of $q$,
which very strongly resembles a definition of local homeomorphism!

\subsection{The reciprocal stability law}

There are several evidences of the fact that the reciprocal stability laws should hold in
an appropriate category $\TT$ of ``topological spaces":
\begin{enumerate}
\item
On one side, the (antiperfect, perfect) factorization in $\Top$ is indeed pullback-stable along local homeomorphisms.
This generalizes the fact that if $P\inc X$ is dense and $O\inc X$ is open, then $P\cap O \inc O$ is dense.
(The related fact that open maps reflect density is taken as a basis for a definition of open maps 
in~\cite{tho04}.) 
\item
For the other stability law, we present three particular cases.
If $X\in\TT$ is a $\rm T_1$ space (that is, its points $x:1\to X$ are in $\M'$), 
the pullback squares below show that the (discrete) fiber $\nh P x$ of the 
etale reflection of a map $p:P\to X$ is given by the components of the fiber space $Px$:
\[  
\xymatrix@R=4pc@C=4pc{
P x  \ar[r]^{e'}\ar[d] &  \nh P x  \ar[r]^{m'}\ar[d]^n  &  1 \ar[d]^x   \\
P \ar[r]^e             &  \nh P \ar[r]^m                &  X             }
\]
Indeed, since $x\in\M'$ and $m\in\M$, also $n\in\M'$ and $m'\in\M$.
Then, by the reciprocal stability, $e\in\E$ implies $e'\in\E$.
Thus the top row gives the discrete reflection of $Px$, that is $\nh P x = \comp Px$.
In fact, giving the quotient topology to the set of fibers components, one gets the etale reflection for
some classes of maps in $\Top$ (see e.g.~\cite{joh}).
\item
If $K\in\TT$ is compact (that is, $K\to 1$ is in $\M'$), the following pullback diagram similarly shows 
that if $\nh P$ is a neighborhood of $P\inc X$ then $K\tm\nh P$ is a neighborhood of $K\tm P \inc K\tm X$: 
\[   
\xymatrix@R=4pc@C=4pc{
K\tm P \ar[r]^{e'}\ar[d] &  K\tm\nh P \ar[r]^{m'}\ar[d]^{n''}  &   K\tm X \ar[r]\ar[d]^{n'}  &  K\ar[d]^n \\
P \ar[r]^e               &  \nh P \ar[r]^m                     &  X   \ar[r]                 &  1          }
\]
In classical terms (in $\Top$), open sets $K\tm O$ form a basis for the open sets in $K\tm X$ containing $K\tm P$,
when $O$ runs through open sets in $X$ containing $P$ (which is of course not true if $K\inc X$ is,
for instance, a straight line in the plane).
\item
Similarly, if $D$ a closed part of $X\in\TT$ (that is, the monomorphism $D\inc X$ is in $\M'$), we have
\[  
\xymatrix@R=4pc@C=4pc{
D\cap P  \ar[r]^{e'}\ar[d] &   D\cap\nh P  \ar[r]^{m'}\ar[d]^{n'}  &  D \ar[d]^n   \\
P \ar[r]^e                 &  \nh P \ar[r]^m                       &  X             }
\]
In classical terms, the open sets $D\cap O$ form a basis for the open sets in $D$ containing $D\cap P$,
when $O$ runs through open sets in $X$ containing $P$ (which is of course not true in general for a non-closed part $D$). 
\item
A consequence of the reciprocal stability law is the exponential law (Proposition~\ref{p2.bfc}):
if $m\in\MX$, $n\in\MoX$ and the exponential $m^n$ exists in $\TX$, then it is in $\MX$
(and conversely).
In particular, if $K$ is compact and $S$ is discrete, $K^S$ is compact and $S^K$ is discrete.
The first one is a consequence, in $\Top$, of Tychonoff theorem, or else it follows from the first point above.
For the second one, note that a compact locally connected space has a finite number of components.
Thus, the compact-open topology shows that $S^K$ is discrete: namely $S^K = S^{\comp K}$.
\item
Again by the exponential law, any ``finite covering" $b$ of $X\in\TT$ (that is, $b\in(\M\cap\M')/X$; see
Section~\ref{hom}) yields a ``$b$-complementation" $\neg_b := b^- :\MX\to\MoX$ (and conversely);
if $\TT$ has an initial object and $0\to 1$ is a ``finite set" in $(\M\cap\M')/1$, we get a complementation  
$\neg_{!_X^*0}$, which generalizes the classical one between open and closed parts in $\Top$. 
Thus the latter is only the trace left on monomorphic parts
of a less perfect but much more pregnant ``duality" between perfect maps and local homemorphisms, 
which is fully expressed by the \rs law.
\end{enumerate}

\subsection{Related work}

``Categorical" or ``universal" topology has a long history and has developed in many different threads
(which in part reflect the variety of the concepts that can be considered as basic in topology itself).
The present work belongs to the one that look for a proper categorical foundation of topology 
via a suitable axiomatization of ``topological" categories
that is, categories $\TT$ whose objects can be effectively considered as topological spaces of kind $\TT$
(in the same sense, say, that the objects of a topos $\T$ can be considered as sets of kind $\T$).

In this direction (but not concerned specifically with classical topology) we have already mentioned the 
fundamental work of Bill Lawvere who develops in several papers an analysis of the objects of a 
category $\C$ by contrasting them with ``discrete" objects; 
furthermore the latter can often be defined inside $\C$ by means of a special object
(for instance, the arrow category in $\Cat$, or a ``tiny" $T\in\TT$ such that $X^T$ is the tangent bundle of $X$).
Here, a similar role may be played by a ``Frechet space" (see Section~\ref{conv}), 
which gives ``discrete" (or ``etale") and ``compact" (or ``perfect") objects at any slice $\TX$.
Anyway, we do not assume exponentiability and the existence of ``interior" right adjoints $\TX\to\MX$ as basic;  
rather, these properties can be considered as further possible axioms (see Section~\ref{top}). 

On the other hand in~\cite{tho} and~\cite{tho04}, it is presented an abstraction of $\Top$ based on closed maps 
and it is developed a great amount of classical topology therein.
In spite of the strictly related basic concepts, however, that approach differs from ours in several respects. 
For instance, we {\em simultaneously} consider perfect maps and local homeomorphisms 
(rather than seeing them as two separated instances of the same abstraction)
and we use \fss to condense their basic properties and reciprocal relationships 
(rather than to handle images of ``subobjects", which are not particularly relevant to us).

Several \fss on $\Top$ have been considered in the literature and 
many of them have been studied in~\cite{joh}, in the context of toposes as generalized spaces. 
Among these, there seems not to be (even in the generalized context)
a natural pair of reciprocally stable factorization systems, so that the question of a 
concrete model for ``balanced topology" remains open.
We mention also the recent work \cite{ane}, concerned with the construction
of a Grothendieck topology associated to a factorization system, especially in the context of
algebraic geometry; there, the etale-proper (or perfect) ``duality"
seems to emerge again in guises related to the present work.

\section{Topological spaces and discrete fibrations}
\label{df}

The analogy between \lhs and \dfs and between \pms and \dofs is one of the main motives of our 
common abtraction of $\Top$ and $\Cat$ as weak bfc's.
We here review two ``explainations" of this analogy.

\subsection{Compactness and discreteness in slices of $\Top$ and of $\Cat$}

Following~\cite{bou}, a space $X\in\Top$ is compact if it is
\begin{enumerate}
\item
quasi-compact, that is all the projections $p:T\tm X\to T$ are closed,
\item
and separated, that is the diagonal map $\Delta:X\to X\tm X$ is closed.
\end{enumerate} 
This definition can be extended to any finitely complete category $\C$ with a functor $(-)\st:\C\to\Top$:
an object $X\in\C$ is compact if it is
\begin{enumerate}
\item
quasi-compact, that is all the maps $p\st :(T\tm X)\st\to T\st$ are closed,
\item
and separated, that is the map $\Delta\st:X\st\to(X\tm X)\st$ is closed.
\end{enumerate} 
With the projection $(-)\st:\TopX\to\Top$, the compact (resp. quasi-compact) objects of $\TopX$
are the perfect (resp. proper) maps to $X$ (see~\cite{bou}).
Replacing ``closed" with ``open" in the above definitions, we similarly get discrete spaces in $\Top$, 
and \lhs (resp. open maps) in $\TopX$.

Considering the functor $(-)\st:\Cat\to\Top$ that sends a category $X$ to the (Alexandroff) space $X\st\in\Top$ of its 
thin reflection, it is easy to see that all categories are quasi-compact, while
the separated, and hence also the compact ones, coincide with the discrete ones.
Composing with $\CatX\to\Cat$, we get a functor $(-)\st:\CatX\to\Top$ giving, as
compact objects, the discrete opfibrations over $X$.
Dually, \lhs in $\CatX$ are the discrete fibrations over $X$.

Of course, by redefining closed parts as monomorphic perfect maps, one gets the ``upward-closed"
full subcategories that is, the closed parts of $X\in\Cat$ are those of $X\st$,
but considered as full subcategory inclusions (and similarly for open parts).

\subsection{Discrete (op)fibrations via orthogonality}

We have just seen a definition of \dofs over $X\in\Cat$ as compact objects in $\CatX$.
But they can be defined more naturally as those functors which are orthogonal to the domain
$s:1\to 2$ of the arrow category; $n:D\to X$ is in $\MoX$ iff any square
\[   
\xymatrix@R=4pc@C=4pc{
1  \ar[r]^a\ar[d]_s            &  D \ar[d]^n \\
2  \ar[r]^l\ar@{..>}[ur]^{l'}  &  X  }
\]
has a unique diagonal.
That is, given an object $a\in D$, any arrow $l$ in $X$ with domain $na$ has a unique lifting (along $n$) 
to an arrow $l'$ with domain $a$.
(Quasi-compact objects are those for which the lifting $l'$ exists but not necessarly unique, as can be checked
by using $T=2$ as test object.)
Dually, \dfs are those functors which are orthogonal to the codomain functor $t:1\to 2$.

\subsection{Perfect maps and \lhs via convergence}
\label{conv}

Perfect maps $n:D\to X$ in $\Top$ can be defined by a similar ``convergence lifting" property (see~\cite{bou});
if $\nu$ is an ultrafilter in $D$ such that $n\nu$ converges to $x\in X$, 
then $\nu$ converges to a unique $a$ over $x$. 

Briefly, \pms are ``ultrafilter opfibrations".
Now, one would expect that (in view of the above ``dual" characterizations)
\lhs $m:O\to X$ in $\Top$ can be dually defined as ``ultrafilter fibrations":
if $a\in O$ and $\xi$ converges to $ma\in X$, 
then there is a unique $\nu$ in $O$ over $\xi$ converging to $a$. 
In fact, in~\cite{hof05} it is shown that \lhs are the {\em pullback stable} ultrafilter fibrations.

We can make the link explicit by assuming that in our ``topological" bfc $\T$
the factorization systems are generated by a ``Frechet" object (instead of the bipointed arrow object of $\Cat$):
\[  
\xymatrix@R=4pc@C=4pc{
1  \ar[r]^e   &  F  &  F' \ar[l]_i   }
\]
(with $e\in\E$ and $i\in\E'$).
$F'$ should be thought of as a ``free sequence", which is included in $F$ as a ``convergent sequence".
In that case, local homeomorphisms and perfect maps {\em are} the Frechet discrete fibrations and opfibrations,
that is maps in $\M$ and $\M'$ are defined by the following unique liftings properties:
\[   
\xymatrix@R=4pc@C=4pc{
1  \ar[r]\ar[d]^e        &  O \ar[d] \\
F  \ar[r]\ar@{..>}[ur]   &  X         }
\qq\qq
\xymatrix@R=4pc@C=4pc{
F'  \ar[r]\ar[d]          &  D \ar[d] \\
F  \ar[r]\ar@{..>}[ur]    &  X         }
\]
(For monomorphic maps, these give the classical convergence characterization of open and closed parts:
$O\inc X$ is open when any ``convergent sequence" in $X$, converging to a point in $O$, is itself (definitively) in $O$; 
$D\inc X$ is closed when for any ``sequence" in $D$, converging to a point $x\in X$, one has $x\in D$.)
\begin{remark}
Convergence is one of the basic ideas of topology.
Its formalization in $\Top$ through ultrafilters has been proved fruitful in several respects, giving often
more intuitive counterparts of definitions and properties.
Beside the above mentioned characterization of \pms and \lhs (and so also of compact spaces, closed and open parts, etc.),
ultrafilters can also be used to define topological spaces themselves and to characterize the exponentiable ones 
(see for instance~\cite{pis3} and~\cite{hof03}).
On the other hand, the use of ultrafilters in topology has some drawbacks;
apart from the lack of constructivity, their practical use is often rather akward 
(as in the proof, in~\cite{hof03}, of the exponentiability of perfect maps).
Furthermore, sometimes the results are not exactly how one could reasonably expect. 
For instance, the fact that an ultrafilter fibrations may not be a local homeomorphism, with the accompanying 
counter-example, appears rather as a flaw of classical topological spaces and of the ultrafilter analysis of convergence,
allowing such ``pathological" spaces.
In our context, it seems to be possible a more direct and intuitive approach to infinitesimal aspects 
and to their analysis via convergence.
\end{remark}

\section{Balanced topology}
\label{top}

Balanced topology is based on the assumption that $\TT$ is a (weak) bfc, whose objects are to be thought
of as some kind of topological spaces, possibly infinitesimal and suitably regular.
We here briefly sketch some properties that follow from this assumption, and hint at some other possible
axioms that may render $\TT$ a better approximation of the idea of a topological category.

\subsection{Terminology and notation}

We refer to the objects of $\TT$ as (topological) {\bf spaces}, to the maps in $\M$ as {\bf \lhs} 
(or also ``discrete" or ``etale" maps) and to maps in $\M'$ as {\bf \pms}.
Maps in $\B = \M\cap\M'$ are the {\bf \fcs}.
Maps in $\E$ (resp. $\E'$) will be called {\bf final} (resp. {\bf initial}) maps 
(altought other names have been used for the latter in $\Top$).
The objects (maps) in $\S := \Mu$ are the {\bf discrete spaces} or {\bf (internal) sets}.
The objects (maps) in $\K := \Mou$ are the {\bf compact spaces}.
The objects (maps) in $\So := \S\cap\K = \B/1$ (the finite coverings of $1$) are the {\bf finite sets}.
(Note that, for $\TT = \Cat$, finite sets may be not... finite.)
The objects (maps) in $\Eu$ are the {\bf connected spaces}. 
Letting $\P X$ be the slice $\TX$ restricted to monomorphisms, $\O X := \P X \cap \MX$ 
are the {\bf open parts} of $X$, and $\D X := \P X \cap \MoX$ are the {\bf closed parts} of $X$. 
The parts in $\D X\cap \O X = \P X \cap \BX$ are {\bf clopen}.
The reflection $\comp:\TT\to\S$ is the {\bf components} functor, and $\comp X$ is the {\bf set of components} of $X$.
A space $X$ is {\bf finite} if its set of components is finite. 
A space $X$ is {\bf separated} if the diagonal $\Delta:X\to X\tm X$ is in $\M'$.
A space is ${\rm\bf T_1}$ if its points are closed.
A space $X$ is {\bf groupoidal} if $\MX = \MoX =\BX$.

If $P\inc X$, its $\EM$-factorization $P\to\nh P\to X$ is the {\bf\ng} of $P$ in $X$.
If it is monomorphic as well, it is both the smallest open part containing $P$ and the biggest 
part of $X$ containing $P$ as a final part (see Corollary~\ref{p3.fs}).

The proposition below simply expresses properties of \fss rephrased in the above language:
\begin{prop}    \label{p1.top}
\begin{itemize}
\item
Perfect maps and \lhs over a space are closed with respect to all the limits which exist in $\TT$; 
in particular finite limits of compact (resp. discrete) spaces are themselves compact (resp. discrete).
\item
If $\TT$ is (finitely) cocomplete, so are \pms and \lhs over a space (in particular, $\K$ and $\S$).
\item
Perfect maps and \lhs ar pullback stable.
The pullback of a \pm (resp. \lh) along a map with a compact (resp. discrete) domain, has itself
a compact (resp. discrete) domain.
(Briefly, a \pm has compact fibers over compact parts.)
\item
Any compact space is separated and ${\rm T_1}$.
Any discrete space has an open diagonal and open points.
\item
The equalizer of two parallel maps to a separated (resp. discrete) space is closed (resp. open). 
\item
Any map between compact spaces is perfect.
Any map between discrete spaces is a local homeomorphism.
\item
A space is connected iff any map to a discrete space is constant.
\item
For any figure $P\to X$ with a connected shape, its \ng $\nh P\to X$ has a connected shape as well;
in particular, any space is locally connected. 
\end{itemize}
\epf
\end{prop}
The following are some ``topologically reasonable" consequences of the \rs laws:
\begin{prop}    \label{p2.top}
\begin{itemize}
\item
Pulling back \ngs along proper maps, one gets \ngs again; in particular, intersecting with closed parts 
or multiplying by compact spaces preserves neighborhoods.
\item
The exponential law holds for exponentiable spaces (see Proposition~\ref{p2.bfc}).
\item
The fiber of a final map over a closed connected part, is connected
(e.g., over points, for ${\rm T_1}$ spaces, or over closures of points if the discrete are separated).
\item
The components of a finite space are connected and clopen.
\end{itemize}
\end{prop}
\pf
Most of these have been already discussed at the beginning of this second part;
for the last one, recall Diagram~(\ref{e7.ten}) and that the points of a
finite (internal) set are clopen.
\epf

\subsection{Further topological axioms}

The following properties hold in $\Top$, and so are possible axioms for $\TT$:
\begin{itemize}
\item
$\TT$ is extensive and $1\in\TT$ is (externally) connected.
\item
$1\in\TT$ is groupoidal and two-valued: $1$ and $0$ are the only (cl)open part of it.
\item
There is a ``Sierpinski" space, which classifies open parts.
\item
There are ``interior" coreflections $\CX\to\MX$, for any $X\in\T$.
\item
Discrete spaces are separated and ${\rm T_1}$.
\item
Perfect maps and \lhs are exponentiable.
\end{itemize}

\subsection{Some homotopical properties}
\label{hom}

Since a perfect local homeomorphism between locally connected topological spaces is a finite covering
(see~\cite{bou}), it is natural to define the class of finite coverings in $\TT$ as $\B = \M\cap\M'$;
then $\BX$ should reflect the \p1-homotopy type of $X\in\TT$.
We say that maps $f:X\to Y$ and $g:Y\to X$ are a {\bf \p1-equivalence} if they induce
an equivalence between $\BX$ and $\BY$.
In particular, a space $I$ is ``simply connected" if it \p1-equivalent to $1\in\T$,
that is, if the finite internal sets inclusion $\S_0\to\BX$ is an equivalence.  

In $\Cat$, we have $\BX\equ\Set^{X'}$, where $X'$ is the groupoidal reflection of $X$.
Thus, for example, an adjunction $f\adj g:X\to Y$ is a \p1-equivalence in $\Cat$,
since it gives an equivalence $f'\adj g':X'\to Y'$.
In particular, a category with a terminal (or initial) object is simply connected.
Another instance of simply connected category is any connected poset.
The following is a consequence of Proposition~\ref{p6.fs}: 
\begin{corol}
The \p1-equivalences have the unique lifting property with respect to finite coverings.
\epf
\end{corol}  
In a topological bfc $\TT$, a map $i:A\inc X$ in $\E$ can be seen as the inclusion of $A$ in one of
its possible neighborhoods in an ``ampler" space (e.g., $X$ itself).
Thus, the following result may be rephrased by saying that an (infinitesimal) \ng of a space $A$ 
which retracts on $A$ has the same \p1-homotopy type of $A$ itself.
\begin{prop}    \label{p1.hom}
A retraction $r,i:A\to X$ with $i\in\E$ is a \p1-equivalence.
\end{prop}
\pf
Let $b:B\to X$ be any finite covering of $X$.
In the diagram below, the left hand square is a pullback and the right hand one is obtained 
by factorizing the map $rb:B\to A$ according to $\EM$:
\[
\xymatrix@R=4pc@C=4pc{
i\st B \ar[r]^{e'}\ar[d]_{i\st b} & B \ar[d]^b\ar[r]^e  & \ex_r B \ar[d]^{b'}  \\
A \ar[r]^i                        & X  \ar[r]^r         & A          }
\]
By the rsl, $e'$ is in $\E$, and so also $e\circ e'$ is in $\E$; 
since $b'\circ e\circ e' = i\st b$ and $i\st b,b'\in\M$, the map $e\circ e'$ is also in $\M$, and so it is an iso.
Thus the adjunction $\ex_r\adj r\st:\M/A\to\MX$ restricts to an adjunction $\ex_r\adj r\st:\B/A\to\BX$.
Since $r\in\E$, again by the rsl the counit $\ex_r r\st b'\to b'$ is an iso for any $b'\in\B/A$. 
It remains to show that the unit $b\to r\st\ex_r b$ is an iso as well.

Pulling back $b'$ along $ri=\id_A$ we get another isomorphism $e'''\circ e'':i\st r\st\ex_r B\to\ex_r B$: 
\[
\xymatrix@R=4pc@C=1pc{
               & i\st r\st\ex_r B\ar[rr]^{e''}\ar[ddl] && r\st\ex_r B\ar[dr]^{e'''}\ar[ddl]^<<<<<<<<<<<{b''}      \\
i\st B \ar[rr]^>>>>>>>>>>{e'}\ar[d]_{i\st b}\ar@{..>}[ur]^s      && B \ar[d]_b\ar[rr]^>>>>>>>>>>e\ar@{..>}[ur]^u  && \ex_r B \ar[d]^{b'}  \\
A \ar[rr]^i                                            && X  \ar[rr]^r                        && A          }
\]
The mediating iso $s$ is easily seen to be a map over $A$ such that $u\circ e' = e''\circ s$,
where $u$ is universally induced to the pullback $r\st\ex_r B$.
Thus the latter is both in $\M$ and in $\E$ that is, it is an isomorphism.
\epf
In particular, any finite covering $b$ of the neighbouring space $X/x$ of a point $x:1\to X$ is ``constant" 
that is, $b \,=\,\,!_{X/x}^*S$, for a finite set $S\in\S_0$:
\begin{corol}  \label{p.hom}
Any space $X\in\TT$ is locally simply connected and any finite covering $b\in\BX$ is ``locally trivial":
pulling back $b$ along a \ng $X/x\to X$ one gets a constant covering.
\epf
\end{corol}

\section{Conclusion of the second part}
\label{con2}

We have shown that assuming that $\TT$ is a (Frechet generated) bfc allows one to capture
several relevant features of topology.
Although this ``version" of topology may appear over-simplified, it has the advantage to offer a direct and intuitive 
approach both to ``local" (or ``infinitely close") aspects of spaces, and also to some ``global" (or homotopical) properties. 
 
In fact, any space $X\in\T$ has a ``left topology" $\MX$ of ``open" figures and a ``right topology" $\MoX$ of ``closed" figures
interacting by the \rs law (wich generalizes the complementation law in classical topology).
Furthermore, $X$ has a ``\p1-homotopy" $\BX$ and also left and right ``cotopologies" $X\bs\E$ and $X\bs\E'$;
for instance, the left cotopology of $1\in\T$ of ``infinitesimal quantities" spaces should be an important
object of study in balanced topology.

Thus, we have sketched a genuinely categorical approach to topology that, we hope, can help
to organize and guide topological thinking and can also offer a topological perspective on category theory.
It remains open the question of what are the proper further axioms for $\TT$,
and if a ``concrete model" of $\TT$ is available.

\begin{refs}

\bibitem[Ad\'amek et al., 2001]{dense} J. Ad\'amek, R. El Bashir, M. Sobral, J. Velebil (2001), 
On Functors which are Lax Epimorphisms, {\em Theory and Appl. Cat.} {\bf 8}, 509-521.

\bibitem[Anel, 2009]{ane} M. Anel (2009), Grothendieck Topologies from Unique Factorization Systems, 
preprint, arXiv:math.AG/0902.1130


\bibitem[Bourbaki, 1961]{bou} N. Bourbaki (1961), {\em Topologie Generale}, Hermann, Paris.

 
 
 

\bibitem[Clementino et al., 1996]{clem96} M.M. Clementino, E. Giuli, W. Tholen (1996), Topology in a Category: Compactness, 
{\em Portugal. Math.} {\bf 53}(4), 397--433.

\bibitem[Clementino et al., 2003]{hof03} M.M. Clementino, D. Hofmann and W. Tholen (2003), 
The Convergence Approach to Exponentiable Maps, {\em Port. Math.} {\bf 60}, 139-160.

\bibitem[Clementino et al., 2004]{tho04} M.M. Clementino, E. Giuli, W. Tholen, (2004), {\em A Functional Approach to Topology}, 
Encyclopedia of Mathematics and its applications, vol. 97, Cambridge University Press, 103-163.

\bibitem[Clementino et al., 2005]{hof05} M.M. Clementino, D. Hofmann and G. Janelidze (2005), 
Local Homeomorphisms via Ultrafilter Convergence, {\em Proc. Amer. Math. Soc.} {\bf 133}, 917-922.

\bibitem[Johnstone, 1982]{joh} P. Johnstone (1982), {\em Factorization Theorems for Geometric Morphisms II}, 
Lecture Notes in Mathematics 915, Springer, Berlin, 216-233.



\bibitem[La Palme et al., 2004]{rey04} M. La Palme Reyes, G.E. Reyes, H. Zolfaghari (2004), 
{\em Generic Figures and their Glueings}, Polimetrica S.a.s., Monza. 

\bibitem[Lawvere, 1966]{law66} F.W. Lawvere (1966), {\em The Category of Categories as a Foundation for Mathematics},
Proceedings of the Conference on Categorical Algebra, La Jolla, 1965, Springer, New York, 1-20.

\bibitem[Lawvere, 1969]{law69} F.W. Lawvere (1969), Adjointness in Foundations,
{\em Dialectica} {\bf 23}, 281-295. Republished in {\em Reprints in Theory and Appl. Cat.}
 
\bibitem[Lawvere, 1970]{law70} F.W. Lawvere (1970), {\em Equality in Hyperdoctrines and the Comprehension Scheme as an Adjoint Functor},
Proceedings of the AMS Symposium on Pure Mathematics, XVII, 1-14. 

\bibitem[Lawvere, 1989]{law89} F.W. Lawvere (1989), {\em Qualitative Distinctions between some Toposes of Generalized Graphs},
Proceedings of the AMS Symposium on Categories in Computer Science and Logic, Contemporary Mathematics, vol. 92, 261-299.

\bibitem[Lawvere, 1992]{law92} F.W. Lawvere (1992), {\em Categories of Space and of Quantities}, The Space of Mathematics,
DeGruyter, Berlin, 14-30.


\bibitem[Lawvere, 2003]{law03} F.W. Lawvere (2003), Foundations and Applications: Axiomatization and Education,
{\em Bull. Symb. Logic} {\bf 9}(2), 213-224.

\bibitem[Lawvere, 2007]{law07} F.W. Lawvere (2007), Axiomatic Cohesion, 
{\em Theory and Appl. Cat.} {\bf 19}, 41-49. 

\bibitem[Maltsiniotis, 2005]{mal} G. Maltsiniotis (2005), Structures d'asph\'ericit\'e, Foncteurs Lisses, et Fibrations, 
{\em Ann. Math. Blaise Pascal} {\bf 12}, 1-39.


\bibitem[Par\'e, 1973]{par} R. Par\'e (1973), Connected Components and Colimits,
{\em J. Pure Appl. Algebra} {\bf 3}, 21-42.

\bibitem[Pisani, 1999]{pis3} C. Pisani (1999), Convergence in Exponentiable Spaces, 
{\em Theory and Appl. Cat.} {\bf 5}, 148-162.

\bibitem[Pisani, 2005]{pis4} C. Pisani (2005), Bipolar Spaces, preprint, arXiv:math.CT/0512194

\bibitem[Pisani, 2007]{pis0} C. Pisani (2007), Components, Complements and the Reflection Formula, 
{\em Theory and Appl. Cat.} {\bf 19}, 19-40. 

\bibitem[Pisani, 2007a]{pis1} C. Pisani (2007a), Components, Complements and Reflection Formulas, preprint, arXiv:math.CT/0701457

\bibitem[Pisani, 2007b]{pis2} C. Pisani (2007b), Categories of Categories, preprint, arXiv:math.CT/0709.0837

\bibitem[Pisani, 2008]{pis} C. Pisani (2008), Balanced Category Theory, 
{\em Theory and Appl. Cat.} {\bf 20}, 85-115. 

\bibitem[Street \& Walters, 1973]{stw} R. Street and R.F.C. Walters (1973), The Comprehensive Factorization of a Functor,
{\em Bull. Amer. Math. Soc.} {\bf 79}(2), 936-941.

\bibitem[Tholen 1999]{tho} W. Tholen (1999), 
A Categorical Guide to Separation Compactness and Perfectness, {\em Homol. Homot. Appl.} {\bf 1}, 147-161.

\end{refs}

\end{document}